\documentstyle[12pt]{article}

\def\R{\bf R}
\def\C{\bf C}
\begin{document}

\begin{center}
{\Large {\bf A Uniformization Theorem Of Complete Noncompact K\" ahler
Surfaces With Positive Bisectional Curvature}} \vskip 5mm Bing--Long Chen*,
Siu--Hung Tang**, Xi--Ping Zhu*

\smallskip 
* Department of Mathematics, Zhongshan University,\\Guangzhou 510275, P. R.
China, and\\ \smallskip 
Institute of Mathematical Sciences, The Chinese University of Hong Kong,
Hong Kong\\\smallskip 
** Institute of Mathematical Sciences, The Chinese University of Hong Kong,
Hong Kong\\ \medskip
\end{center}

\baselineskip=20pt

\begin{abstract}
In this paper, by combining techniques from Ricci flow and algebraic
geometry, we prove the following generalization of the classical
uniformization theorem of Riemann surfaces. Given a complete noncompact 
complex two dimensional K\"ahler manifold $M$ of positive and bounded holomorphic
bisectional curvature, suppose its geodesic balls have Euclidean
volume growth and its scalar curvature decays to zero at infinity in the 
average sense, then $M$ is biholomorphic to $\C^2$. During the 
proof, we also discover an interesting gap phenomenon which says 
that a K\"ahler manifold as above automatically has quadratic
curvature decay at infinity in the average sense.
\end{abstract}

\vskip 0.5mm

\section*{\S1. Introduction}

\setcounter{section}{1} \setcounter{equation}{0}

\qquad One of the most beautiful results in complex analysis of
one variable is the
classical uniformization theorem of Riemann surfaces which
states that a simply connected Riemann surface is biholomorphic 
to either the
Riemann sphere, the complex line or the open unit disc.
Unfortunately, a direct analog of this beautiful result 
to higher dimensions does not exist. For example,
there is a vast variety of
biholomorphically distinct complex structures on ${\R}^{2n}$ 
for $n>1$, a fact 
which was already known to Poincar\'e (see \cite{BSW},
\cite{Fe} for a modern treatment). 
Thus, in order to characterize 
the standard complex structures for
higher dimensional complex manifolds, one must impose more
restrictions on the manifolds.

From the point of view of differential geometry,
one consequence of the uniformization theorem is that a positively
curved compact or noncompact Riemann surface must be biholomorphic
to the Riemann sphere or the complex line respectively. 
It is thus natural to ask whether there is similar    
characterization for higher dimensional complete K\"ahler 
manifold with positive "curvature". 
%Note that one of the consequences of the classical uniformization theorem 
%is the uniqueness of the complex structure of a positive curved
%Riemann surface. 
That such a characterization exists in the case of 
compact K\"ahler manifold  
is the famous Frankel conjecture which says that
a compact K\"ahler manifold of positive holomorphic bisectional 
curvature is biholomorphic to a
complex projective space. This conjecture was solved by
Andreotti--Frankel \cite{Fra} and Mabuchi \cite{Mab} in
complex dimensions two and three respectively and the general 
case was then solved by Mori \cite{Mor}, and
Siu--Yau \cite{SiY} independently.
In this paper, we are thus interested in complete noncompact K\"ahler 
manifolds with positive holomorphic bisectional curvature. The 
following conjecture provides the main impetus.
%the following conjecture concerning the complex
%structure of complete noncompact K\"ahler manifolds with 
%positive bisectional curvature
%is our main interest.
%provide such a natural higher--dimensional analogue.
%\vskip 3mm {\bf 
%\underline{Conjecture I}} (Frankel\cite{Fra})\quad 
%A compact K\"ahler
%manifold of positive holomorphic bisectional 
%curvature is biholomorphic to a
%complex projective space.
\vskip 3mm
{\bf \underline{Conjecture}}
(Green--Wu \cite{GW2}, Siu\cite{Si}, Yau \cite{Y2})\quad 
A complete noncompact\\K\"ahler manifold of
positive holomorphic bisectional curvature is biholomorphic 
to a complex Euclidean space.
\vskip 3mm 
%The first conjecture is known as Frankel
%conjecture. Andreotti--Frankel \cite{Fra} proved that Frankel conjecture is
%true when the K\"ahler manifold has complex dimension two. The case of
%complex dimension three was solved by Mabuchi \cite{Mab}. And the general
%dimension case of Frankel conjecture was proved by Mori \cite{Mor}, and
%Siu--Yau \cite{SiY}. This says, Conjecture I has been completely solved.

In contrary to the compact case, very little is known 
about this conjecture. The first result in this direction
is the following isometrically embedding theorem.

\vskip 3mm{\bf \underline{Theorem}} (Mok--Siu--Yau\cite{MSY}, Mok\cite{Mo1}) 
\quad 
Let $M$ be a complete noncompact K\"ahler manifold of nonnegative
holomorphic bisectional curvature of complex dimension $n\geq 2$. 
Suppose there exist positive constants $C_1$, $C_2$ such that for 
a fixed base point $x_0$ and some $\varepsilon>0$,
%\vskip 0.5mm 
$$
\begin{array}{lll}
\bigbreak 
\mbox{(i)} \qquad & \mbox{Vol}\,(B(x_0,r))\geq C_1r^{2n}\  
& \qquad \qquad 0\leq r<+\infty \ ,\qquad \qquad \\ 
\mbox{(ii)} \qquad & R(x)\leq \displaystyle \frac{C_2}
{1+d^{2+\varepsilon}(x_0,x)}\  & \qquad
\qquad \mbox{on}\quad M\ ,\qquad \qquad 
\end{array}
$$
where $\mbox{Vol}\,(B(x_0,r))$ denotes the volume of the 
geodesic ball $B(x_0,r)$ centered at $x_0$ with radius $r$, 
$R(x)$ denotes the scalar curvature and $%
d(x_0,x)$ denotes the geodesic distance between $x_0$ and $x$. 
Then $M$ is isometrically biholomorphic to $\C^n$ 
with the flat metric.

\vskip 3mm 
Their method is to consider the Poincar\'e--Lelong equation 
$\sqrt{-1}\partial \overline{\partial }u= \mbox{Ric}$.
Under the condition (ii) that the curvature has faster than 
quadratic decay, they proved the existence of a bounded solution 
$u$ to the Poincar\'e--Lelong equation. By virtue of Yau's 
Liouville theorem on complete manifolds with nonnegative Ricci 
curvature, this bounded plurisubharmonic function $u$ must
be constant and hence the Ricci curvature must be identically zero.
This implies that the K\"ahler metric is flat because of the nonnegativity
of the holomorphic bisectional curvature. However, this argument
breaks down if the faster than quadratic decay condition (ii) is
weaken to a quadratic decay condition. In this case, although we
can still solve the Poincar\'e--Lelong equation with 
logarithmic growth, the boundedness of the solution
can no longer be guaranteed.

In \cite{Mo1}, Mok also developed a 
general scheme for compactifying complete
K\"ahler manifolds of positive holomorphic bisectional curvature. This
allowed him to
obtain the following improvement of the above theorem.
\vskip 3mm
{\bf {\underline{Theorem}}} (Mok\cite{Mo1})\quad 
Let $M$ be a complete noncompact
K\"ahler manifold of complex dimension $n$ with 
positive holomorphic bisectional curvature. Suppose there 
exist positive constants $C_1$, $C_2$
such that for a fixed base point $x_0$,
%\vskip 0.5mm 
$$
\begin{array}{lll}
\bigbreak
\mbox{(i)} \qquad & \mbox{Vol}\,(B(x_0,r))\geq C_1r^{2n}\  
& \qquad \qquad 0\leq r<+\infty \ ,\qquad \qquad \\ 
\mbox{(ii)}^{\prime} \qquad & 0<R(x)\leq \displaystyle 
\frac{C_2}{1+d^2(x_0,x)}\ 
& \qquad \qquad \mbox{on} \quad M\ ,\qquad \qquad 
\end{array}
$$
then $M$ is biholomorphic to an affine algebraic variety. 
Moreover, if in addition the complex dimension $n=2$ and 

(iii)\qquad the Riemannian sectional curvature of $M$ is positive,

\noindent then $M$ is biholomorphic to $\C^2$.

\vskip 3mm 
To the best of our knowledge, the above result
of Mok and its slight improvements by To \cite{T}, 
and Chen-Zhu \cite{CZ2} 
are the best results in complex dimension two of the 
above stated conjecture.
%in complex two--dimensional
%manifolds are the only partial answers to Conjecture II. 
Here, we would also like to recall the remark pointed out in
\cite{CZ2} that there is a gap in the proof of Shi \cite{Sh3}
(see \cite{CZ2} for more explanation) which
would otherwise constitute a better result than that of Mok \cite{Mo1}.

%We remark that
%Shi announced in \cite{Sh2} that under the same assumptions as in the above
%Mok's theorem, the K\"ahler manifold $M$ is biholomorphic to $\C^n$. This is
%actually the main theorem of his Ph. D. thesis \cite{Sh3}. As pointed out in 
%\cite{CZ2}, there is gap in the proof \cite{Sh3}. In \cite{Sh3}, Shi
%constructed a flat metric metric on $M$. But the flat metric may be
%incomplete. Thus one can only get the biholomorphic embedding of $M$ into a
%domain of $\C^n$. Since the set of complex structures on a domain of $\R%
%^{2n} $ $(n>1)$ is huge, this result of Shi \cite{Sh3} is far away to give
%an affirmative answer to Conjecture II.

%Now, our principal result in this paper is the following 
%improvement of the result of Mok in the case of
%complex dimension two.
In this paper, we consider only the case of complex dimension
two. Our principal result is the following
\vskip 3mm {\bf \underline{Main Theorem}} \quad 
Let $M$ be a complete noncompact complex two--dimen-
sional K\"ahler manifold of positive and bounded holomorphic bisectional
curvature. Suppose there exists a positive constant $C_1$ such that for a
fixed base point $x_0$, we have
%\vskip 0.5mm 
$$
\begin{array}{ll}
\bigbreak 
\mbox{(i)}\qquad & \mbox{Vol}(B(x_0,r))\geq C_1r^4\ 
\qquad \ 0\leq r<+\infty \ , \\ 
\mbox{(ii)}^{\prime \prime }\qquad & \displaystyle \lim 
\limits_{r\rightarrow + \infty }\frac 1{\mbox{Vol}(B(x_0,r))}\int_{B(x_0,r)}R(x)dx=0\ , 
\end{array}
$$
then $M$ is biholomorphic to $\C^2$.
\vskip 3mm 
We remark that the condition 
$\mbox{(ii)}^{\prime \prime}$ means that the scalar curvature tends to 
zero at infinity in average sense. In view of the classical 
Bonnet--Myers theorem, this condition is almost necessary to 
make sure that the manifold is noncompact
under our positive bisectional curvature condition.
It is also clear that the pointwise decay condition 
$$
\lim \limits_{d(x,x_0)\rightarrow + \infty }R(x)=0
$$ 
is stronger than $\mbox{(ii)}^{\prime \prime }.$

The proof of the Main Theorem will be divided into three parts. 
In the first part, we will show that $M$ is a Stein manifold
homeomorphic to $\R^4.$ For this, we evolve the K\"ahler metric 
on $M$ by the Ricci flow first studied by Hamilton. Note that
the underlying complex structure of $M$
is unchanged under the Ricci flow, 
thus we can replace the K\"ahler metric in our main
theorem by any one of the evolving metric. The advantage
is that, in our case, the properties of the 
evolving metric are improving during the flow. 
Moreover,
we know that the
Euclidean volume growth condition (i) as well as the positive 
holomorphic bisectional curvature condition are perserved by the evolving 
metric. More importantly, by a blow up and blow down argument as
in \cite{CZ1}, we can prove that the curvature of the evolving metric decays
linearly in time. This implies that the injectivity radius of 
the evolving metric is getting bigger and bigger and any geodesic ball 
with radius less than half of the injectivity radius is almost 
pseudoconvex. By a perturbation argument as in \cite{CZ2}, we 
are then able to modify these geodesic balls to
a sequence of exhausting pseudoconvex domains of $M$ such that any
two of them form a Runge pair. From this, it follow readily 
that $M$ is a Stein manifold homeomorphic to $\R^4$.

In the second part of the proof, we consider the algebra 
$P(M)$ of holomorphic functions of polynomial growth on $M$
and we will prove that its quotient field has transcendental degree 
two over $\C$. For this, we first need to construct
two algebraically independent holomorphic functions 
in the algebra $P(M)$. Using the $L^2$ estimates of 
%$\overline{\partial }$--operator of 
Andreotti--Vesentini \cite{AV} and
H\"omander \cite{Ho}, it suffices to construct a strictly
plurisubharmonic function of logarithmic growth on $M$. Now, 
if the scalar curvature decays in space at least quadratically, 
it was known from \cite{MSY}, \cite
{Mo1} that such a strictly plurisubharmonic function of logarithmic growth
can be obtained by solving the Poincar\'e--Lelong equation, as we mentioned 
before. However, our decay
assumption $\mbox{(ii)}^{\prime \prime }$ is too weak to 
apply their result directly. 
To resolve this difficulty, we make use of the Ricci flow to
verify a new gap phenomenon which was already predicted by Yau 
in \cite{Y3}. More
explicitly, by using the time decay estimate of evolving metric in the
previous part, we prove that the curvature of the initial
metric must decay quadratically in space in certain average sense. 
Fortunately, this turns out to be enough to insure the existence of 
a strictly plurisubharmonic function of logarithmic growth. Next, 
by using the time decay estimate and the
injectivity radius estimate of the evolving metric, we prove that
the dimension of the space of holomorphic functions
in $P(M)$ of degree at most $p$ is bounded by a constant times $p^2$. 
Combining this with the existence of two algebraically independent 
holomorphic functions in $P(M)$ as above, we can prove that  
the quotient field $R(M)$ of $P(M)$ has transcendental degree 
two over $\C$ by a classical argument of Poincar\'e--Siegel.
%and a classical argument of Poincar\'e--Siegel, we can prove
%immediately that the quotient field $R(M)$ of $P(M)$ is a finite extension
In other words, $R(M)$ is a finite extension
field of some ${\C}(f_1,f_2)$, where $f_1$, $f_2 \in P(M)$ are algebraically
independent over $\C$. Then, from the primitive element theorem, we have $R(M)={\C}
(f_1,f_2,g/h)$ for some $g$, $h\in P(M)$. Hence the mapping $F:M\rightarrow {\C
}^4$ given by $F=(f_1,f_2,g,h)$ defines, in an appropriate sense, a
birational equivalence between $M$ and some irreducible affine algebraic
subvariety $Z$ of $\C^4$.

In the last part of proof, we will basically follow the approach of Mok in 
\cite{Mo1} and \cite{Mo3} to establish a biholomorphic map from $M$ onto a
quasi--affine algebraic variety by desingularizing the map $F$. Our
essential contribution in this part is to establish uniform estimates on the
multiplicity and the number of irreducible components of the zero divisor of
a holomorphic function in $P(M)$. Again, the time decay estimate of the Ricci
flow plays a crucial role in the arguments. Based on these estimates, we
can show that the mapping $F:M\rightarrow Z$ is almost surjective in the
sense that it can miss only a finite number of subvarieties in $Z$, and can be
desingularized by adjoining a finite number of holomorphic functions of
polynomial growth. This completes the proof that $M$ is a quasi--affine 
algebraic variety. Finally, by combining with the fact that $M$ is 
homeomorphic to $\R^4$,
we conclude that $M$ is indeed biholomorphic to $\C^2$ by
a theorem of Ramanujam \cite{R} on algebraic surfaces.

This paper contains eight sections. From Sections 2 to 4, we study the Ricci
flow and obtain several geometric estimates for the evolving metric. In Section
5, we show that the two dimensional K\"ahler manifold is homeomorphic to $\R
^4$ and is a Stein manifold. Based on the estimates on the Ricci flow, a
space decay estimate on the curvature and the existence of a strictly
plurisubharmonic function of logarithmic growth are obtained in Section 6.
In Section 7, we establish uniform estimates on the multiplicity and
the number of irreducible components of the zero divisor of a holomorphic
function of polynomial growth. Finally, in Section 8 we construct a
biholomorphic map from the K\"ahler manifold onto a quasi--affine algebraic
variety and complete the proof of the Main Theorem.

We are grateful to Professor L. F. Tam for many helpful discussions and
Professor S. T. Yau for his interest and encouragement.

\section*{\S2. Preserving the volume growth}

\setcounter{section}{2} \setcounter{equation}{0}

\qquad Let $(M,g_{\alpha \overline{\beta }})$ be a complete, 
noncompact K\"ahler surface (i.e., a K\"ahler manifold of 
complex dimension two) satisfying all the assumptions in
the Main Theorem. 
%with bounded and positive
%holomorphic bisectional curvature. Suppose 
%$(M,g_{\alpha \overline{\beta }})$
%satisfies the following conditions:
%\vskip 0.5mm 
%$$
%\begin{array}{ll}
%\bigbreak (i)\qquad & Vol(B(x_0,r))\geq C_1r^4\ ,\qquad \ 0\leq r<+\infty \
%, \\ 
%(ii)^{\prime \prime }\qquad & \displaystyle \lim \limits_{r\rightarrow
%+\infty }\frac 1{Vol(B(x_0,r))}\int_{B(x_0,r)}R(x)dx=0\ , 
%\end{array}
%$$
%where $x_0$ is some fixed point in $M$, $C_1>0$, $B(x_0,r)$ is the geodesic
%ball of radius $r$ with center at $x_0$, and $R(x)$ is the scalar curvature.
We evolve the metric $g_{\alpha \overline{\beta }}$ according to the
following Ricci flow equation 
\begin{equation}
\label{2.1}\left\{ 
\begin{array}{ll}
\bigbreak \displaystyle \frac{\partial g_{\alpha \overline{\beta }}(x,t)}{%
\partial t}=-R_{\alpha \overline{\beta }}(x,t)\  & \qquad x\in M\ \quad
t>0\ , \\ 
g_{\alpha \overline{\beta }}(x,0)=g_{\alpha \overline{\beta }}(x)\  & 
\qquad x\in M\ ,
\end{array}
\right. 
\end{equation}
where $R_{\alpha \overline{\beta }}(x,t)$ denotes the Ricci curvature tensor
of the metric $g_{\alpha \overline{\beta }}(x,t)$.

Since the curvature of the initial metric is bounded, it is known from \cite
{Sh1} that there exists some $T_{\max }>0$ such that (\ref{2.1}) has a
maximal solution on $M\times [0,T_{\max })$ with either $T_{\max
}=+\infty $ or the curvature becomes unbounded as $t\rightarrow T_{\max }$
when $T_{\max }<+\infty $. By using the maximum principle, one knows ( see
Mok \cite{Mo2}, Hamilton \cite{Ha1}, or Shi \cite{Sh4} ) that the positivity
of holomorphic bisectional curvature and the K\"ahlerity of 
$g_{\alpha \overline{\beta }}$ are preserved under
the evolution of (\ref{2.1}). In particular, the Ricci curvature remains
positive.

Our first result for the solution of the Ricci flow 
(\ref{2.1}) is the following proposition.
\vskip 3mm{\bf \underline{Proposition 2.1}} \quad
Suppose $(M,g_{\alpha \overline{\beta }})$ is assumed as above. 
Then the maximal volume growth condition (i) is preserved 
under the evolution of (\ref{2.1}), i.e., 
\begin{equation}
\label{2.2}
\mbox{Vol}_t(B_t(x,r))\geq C_1r^4\ \qquad \mbox{for all} \quad r>0\ ,\quad
x\in M\ , 
\end{equation}
with the same constant $C_1$ as in condition (i). Here, 
$B_t(x,r)$ is the geodesic ball of radius $r$ with 
center at $x$ with respect to the metric 
$g_{\alpha \overline{\beta }}(\cdot ,t)$, and the
volume $\mbox{Vol}_t$ is also taken with respect to the metric 
$g_{\alpha \overline{\beta }}(\cdot ,t)$.

\vskip 3mm{\bf \underline{Proof.}} \quad 
Define a
function $F(x,t)$ on $M\times [0,T_{\max })$ as follows,
$$
F(x,t)=\log \frac{\det \left( g_{\alpha \overline{\beta }}(x,t)\right) }{%
\det \left( g_{\alpha \overline{\beta }}(x,0)\right) }\ . 
$$
By (\ref{2.1}), we have
\begin{eqnarray}
\label{2.3}
\frac{\partial F(x,t)}{\partial t} & = & 
g^{\alpha \overline{\beta }}(x,t)\cdot \frac \partial 
{\partial t}g_{\alpha \overline{\beta }}(x,t)\nonumber \\
& = &-R(x,t) \leq 0 \ ,
\end{eqnarray}
which implies that $F(\cdot ,t)$ is nonincreasing in time.
Since $R_{\alpha \overline{\beta }}(x,t)\geq 0$, 
we know from (\ref{2.1})
that the metric is shrinking in time. In particular, 
\begin{equation}
\label{2.4}
g_{\alpha \overline{\beta }}(x,t)\leq g_{\alpha \overline{\beta }}(x,0)\ 
\qquad \mbox{on} \quad M\times [0,T_{\max })\ . 
\end{equation}
This implies that
\begin{eqnarray}
\label{2.5}
e^{F(x,t)}R(x,t)&=&g^{\alpha \overline{\beta }}(x,t)
R_{\alpha \overline{\beta }}(x,t)\cdot 
\frac{\det \left( g_{\gamma \overline{\delta }}(x,t)\right) }
{\det \left( g_{\gamma \overline{\delta }}(x,0)\right) } \nonumber \\
& \leq & g^{\alpha \overline{\beta }}(x,0)R_{\alpha \overline{\beta }}(x,t) \nonumber \\
&   =  & g^{\alpha \overline{\beta }}(x,0)\left( R_{\alpha \overline{\beta }}(x,t)-R_{\alpha \overline{\beta }}(x,0)\right) +R(x,0)\nonumber\\
&   =  & -\bigtriangleup _0F(x,t)+R(x,0)\ ,
\end{eqnarray}
where $\bigtriangleup _0$ denotes the Laplace operator with respect to the
initial metric $g_{\alpha \overline{\beta }}(x,0)$ and $R(x,t)$ denotes the
scalar curvature of the metric $g_{\alpha \overline{\beta }}(x,t)$.

Combining (\ref{2.3}) and (\ref{2.5}) gives 
\begin{equation}
\label{2.6}e^{F(x,t)}\frac{\partial F(x,t)}{\partial t}\geq \bigtriangleup
_0F(x,t)-R(x,0)\ \qquad \mbox{on} \quad M\times [0,T_{\max })\ . 
\end{equation}

Next, we introduce a cutoff function which will be used several times in 
this paper. Now, as the Ricci curvature of the
initial metric is positive, we know from Schoen and Yau (Theorem 1.4.2 in 
\cite{ScY}) or Shi \cite{Sh4} that there exists a positive 
constant $C_3$ depending only on the dimension such that for 
any fixed point $x_0\in M$ and
any number $0<r<+\infty $, there exists a smooth function 
$\varphi (x)$ on $M $ satisfying 
\begin{equation}
\label{2.7}\left\{ 
\begin{array}{l}
\bigbreak\displaystyle  e^{-C_3\left( 1+\frac{d_0(x,x_0)}r\right) }\leq
\varphi (x)\leq e^{-\left( 1+\frac{d_0(x,x_0)}r\right) }\ , \\ 
\bigbreak      \displaystyle \left| \nabla \varphi \right| _0(x)\leq \frac{%
C_3}r\varphi (x)\ , \\ \displaystyle \left| \bigtriangleup _0\varphi \right|
(x)\leq \frac{C_3}{r^2}\varphi (x)\ , 
\end{array}
\right. 
\end{equation}
for all $x\in M$, where $d_0(x,x_0)$ is the distance between $x$ and $x_0$
with respect to the initial metric $g_{\alpha \overline{\beta }}(x,0)$ and $%
\left| \cdot \right| _0$ stands for the corresponding $C^0$ norm of the
initial metric $g_{\alpha \overline{\beta }}(x,0)$.

Combining (\ref{2.6}) and (\ref{2.7}), we obtain
\begin{eqnarray}
\frac \partial {\partial t}\int_M\varphi (x)e^{F(x,t)}dV_0&\geq&\int_M\left( \bigtriangleup _0F(x,t)-R(x,0)\right) \varphi (x)dV_0\nonumber\\
&\geq&\frac{C_3}{r^2}\int_MF(x,t)\varphi (x)dV_0-\int_MR(x,0)\varphi (x)dV_0\ ,\nonumber
\end{eqnarray}
where $dV_0$ denotes the volume element of the initial metric $g_{\alpha 
\overline{\beta }}(x,0)$.

Recall that $F(\cdot ,t)$ is nonincreasing in time and $F(\cdot ,0)\equiv 0$%
. We integrate the above inequality from $0$ to $t$ to get 
\begin{equation}
\label{2.8}\int_M\varphi (x)\left( 1-e^{F(x,t)}\right) dV_0\leq \frac{C_3t}{%
r^2}\int_M\left( -F(x,t)\right) \varphi (x)dV_0+t\int_MR(x,0)\varphi
(x)dV_0\ . 
\end{equation}
Since the metric is shrinking under the Ricci flow, we have
$$
B_t(x_0,r)\supset B_0(x_0,r)\ \qquad \mbox{for} \quad t\geq 0\ ,
\quad 0<r<+\infty \ , 
$$
and
\begin{eqnarray}
\label{2.9}
\mbox{Vol}_t(B_t(x_0,r))& \geq & \mbox{Vol}_t(B_0(x_0,r)) \nonumber \\
& = & \int_{B_0(x_0,r)}e^{F(x,t)}dV_0 \nonumber \\
& = & \mbox{Vol}_0(B_0(x_0,r))+\int_{B_0(x_0,r)}\left( e^{F(x,t)}-1\right) dV_0\ .
\end{eqnarray}
Then by (\ref{2.7}) and (\ref{2.8}), the last term in (\ref{2.9}) satisfies%
\begin{eqnarray}
\label{2.10}
\int_{B_0(x_0,r)}\left( e^{F(x,t)}-1\right) dV_0&\geq&e^{2C_3}\int_M\left( e^{F(x,t)}-1\right) \varphi (x)dV_0\nonumber\\
&\geq&\frac{C_3e^{2C_3}t}{r^2}\int_MF(x,t)\varphi (x)dV_0\nonumber\\
&&-e^{2C_3}t\int_MR(x,0)\varphi (x)dV_0\ .
\end{eqnarray}
To estimate the two terms of the right hand side of (\ref{2.10}), 
we consider any fixed $%
T_0<T_{\max }$. Since the curvature is uniformly bounded on $M\times [0,T_0]$%
, it is clear from the equation (\ref{2.3}) that $F(x,t)$ is also uniformly
bounded on $M\times [0,T_0]$.

Set
$$
A=\sup \left\{ \left. \left| F(x,t)\right| \right| x\in M\ ,\ t 
\in[0,T_0]\right\}  
$$ 
and
$$
M(r)=\sup \limits_{a\geq r}\frac
1{\mbox{Vol}_0\left( B_0(x_0,a)\right) }\int_{B_0(x_0,a)}R(x,0)dV_0\ . 
$$
Then condition $\mbox{(ii)}^{\prime \prime }$ says that $M(r)\rightarrow 0$ as $%
r\rightarrow +\infty $. By using the standard volume comparison theorem and (%
\ref{2.7}), we have
\begin{eqnarray}
\label{2.11}
\int_MR(x,0)\varphi (x)dV_0&\leq&\int_MR(x,0)e^{-\left( 
1+\frac{d_0(x,x_0)}r\right) }dV_0 \nonumber \\
& = & \int_{B_0(x_0,r)}R(x,0)e^{-\left( 
1+\frac{d_0(x,x_0)}r\right) }dV_0 \nonumber \\
&    & + \sum\limits_{k=0}^\infty \int_{B_0(x_0,2^{k+1}r)
\backslash B_0(x_0,2^kr)}R(x,0)e^{-\left( 
1+\frac{d_0(x,x_0)}r\right) }dV_0 \nonumber \\
& \leq & \int_{B_0(x_0,r)}R(x,0)dV_0 + \sum\limits_{k=0}^\infty 
e^{-2^k}\left( 2^{k+1}\right) ^4 \cdot \nonumber \\
&   & \frac{\mbox{Vol}_0\left( B_0(x_0,r)\right) }{\mbox{Vol}_0
\left( B_0(x_0,2^{k+1}r)\right) }\int_{B_0(x_0,2^{k+1}r)}R(x,0)dV_0\nonumber\\
& \leq & C_4\cdot M(r)\cdot \mbox{Vol}_0\left( B_0(x_0,r)\right) \ ,
\end{eqnarray}
and similarly
\begin{eqnarray}
\label{2.12}
\int_M\varphi (x)dV_0&\leq&\int_{B_0(x_0,r)}e^{-\left( 1+\frac{d_0(x,x_0)}r\right) }dV_0\nonumber\\
&&+\sum\limits_{k=0}^\infty \int_{B_0(x_0,2^{k+1}r)\backslash B_0(x_0,2^kr)}e^{-\left( 1+\frac{d_0(x,x_0)}r\right) }dV_0\nonumber\\
&\leq&\mbox{Vol}_0\left( B_0(x_0,r)\right) +\sum\limits_{k=0}^\infty e^{-2^k}\left( 2^{k+1}\right) ^4\cdot \mbox{Vol}_0\left( B_0(x_0,r)\right)\nonumber\\
&\leq&C_4\mbox{Vol}_0\left( B_0(x_0,r)\right) \ ,
\end{eqnarray}
where $C_4$ is some positive constant independent of $r$.

Substituting (\ref{2.10}), (\ref{2.11}) and (\ref{2.12}) into (\ref{2.9})
and dividing by $r^4$, we obtain%
\begin{eqnarray}
\frac{\mbox{Vol}_t(B_t(x_0,r))}{r^4}&\geq&\frac{\mbox{Vol}_0(B_0(x_0,r))}{r^4}-\frac{C_3e^{2C_3}AT_0}{r^2}\left( C_4\frac{\mbox{Vol}_0(B_0(x_0,r))}{r^4}\right)\nonumber\\
&& -e^{2C_3}T_0\left( C_4M(r)\cdot \frac{\mbox{Vol}_0(B_0(x_0,r))}{r^4}\right)\nonumber\\
&\geq&C_1-\frac{C_3e^{2C_3}AT_0\cdot C_4C_1}{r^2}-e^{2C_3}T_0C_4C_1\cdot M(r)\nonumber
\end{eqnarray}
by condition (i). 
Then letting $r\rightarrow +\infty $,
we deduce that
$$
\lim \limits_{r\rightarrow +\infty }\frac{\mbox{Vol}_t\left( B_t(x_0,r)\right) }{r^4%
}\geq C_1\ . 
$$
Hence, by using the standard volume comparison theorem we have
$$
\mbox{Vol}_t\left( B_t(x,r)\right) \geq C_1r^4\ \quad \mbox{for all}\ x\in M\ ,\ 0\leq
r<+\infty,  \,\, t\in [0,T_0]\ . 
$$
Finally, since $T_0<T_{\max }$ is arbitrary, this completes the proof of the
proposition.
\hfill $\Box$

\section*{\S3. Singularity Models}

\setcounter{section}{3} \setcounter{equation}{0}

\qquad In this and the next section we will continue our study of the
Ricci flow (\ref{2.1}). We will use rescaling arguments to analyse the
behavior of the solution of (\ref{2.1}) near the maximal time $T_{\max }.$

First of all, let us recall some basic terminologies. 
According to Hamilton ( see for example, definition 16.3 in \cite{Ha3}),
a solution to the
Ricci flow, where either the manifold is compact or at each time $t$ the
evolving metric is complete and has bounded curvature, is called a
singularity model if it is not flat and is 
of one of the following three types. Here, we have used $Rm$ to denote 
the Riemannian curvature tensor.

\begin{description}
\item[Type I:]  The solution exists for $-\infty <t<\Omega $ for some $
0<\Omega <+\infty $ and 
$$
|Rm|\leq \frac \Omega {\Omega -t} 
$$
\hspace{6mm} everywhere with equality somewhere at $t=0$;

\item[Type II:]  The solution exists for $-\infty <t<+\infty $ and 
$$
|Rm|\leq 1 
$$
\hspace{8mm} everywhere with equality somewhere at $t=0$;

\item[Type III:]  The solution exists for $-A<t<+\infty $ for some 
$0<A<+\infty$ and 
$$
|Rm|\leq \frac A{A+t} 
$$
\hspace{10mm} everywhere with equality somewhere at $t=0$.
\end{description}

The singularity models of Type I and II are called ancient 
solutions in the sense that the existence time interval of 
the solution contains $(-\infty ,0]$.

Next, we recall the local injectivity radius estimate 
of Cheeger, Gromov and Taylor \cite{CGT}. Let $N$ be a 
complete Riemannian manifold of dimension $m$ with 
$\lambda \leq \mbox{sectional curvature of}\, N \leq \Lambda $ 
and let $r$ be a positive constant satisfying 
$r\leq \frac \pi {4\sqrt{\Lambda }}$ if $\Lambda >0$, 
then the injectivity radius of $N$ at a point $x$
is bounded from below as follows,
$$
\mbox{inj}_N(x) \geq r\cdot \frac{\mbox{Vol}(B(x,r))}{\mbox{Vol}
(B(x,r))+V^m_{\lambda}(2r)}\ , 
$$
where $V^m_{\lambda}(2r)$ denotes the volume of a ball 
with radius $2r$ in the $m$ dimensional model space $V^m_{\lambda}$ 
with constant sectional curvature $\lambda .$
In particular, it implies that for a complete Riemannian 
manifold $N$ of dimension $4$ with sectional curvature bounded 
between $-1$ and $1$, the injectivity radius at a point $x$ 
can be estimated as  
\begin{equation}
\label{3.1}
\mbox{inj}_N(x) \geq \frac 12\cdot \frac{\mbox{Vol}(B(x,\frac 12))}
{\mbox{Vol}(B(x,\frac 12))+V} 
\end{equation}
for some absolute positive constant $V$. Furthermore, if in addition $N$
satisfies the maximal volume growth condition 
$$
\mbox{Vol} \left( B(x,r)\right) \geq C_1r^4\ ,\quad 0\leq r<+\infty \ , 
$$
then (\ref{3.1}) gives 
\begin{equation}
\label{3.2}
\mbox{inj}_N(x) \geq \beta >0 
\end{equation}
for some positive constant $\beta $ depending only on $C_1$ and $V.$

Now, return to our setting. 
Let $(M,g_{\alpha \overline{\beta }})$ be a complete, 
noncompact K\"ahler surface satisfying the same assumptions 
as in the Main Theorem and let $g_{\alpha \overline{\beta }
}(x,t)$ be the solution of the Ricci flow (\ref{2.1}) on 
$M\times [0,T_{\max})$. Denote
$$
R_{\max }(t)=\sup \limits_{x\in M}R(x,t)\ . 
$$
We have shown in Proposition 2.1 that the solution 
$g_{\alpha \overline{ \beta }}(\cdot ,t)$ satisfies the 
same maximal volume growth condition (i)
as the initial metric. Since condition (i) is invariant under
rescaling of metrics, by a simple rescaling argument, 
we get the following injectivity radius estimate for the 
solution $g_{\alpha \overline{\beta }}(\cdot ,t),$
\begin{equation}
\label{3.3}
\mbox{inj}(M,g_{\alpha \overline{\beta }}(\cdot ,t))\geq \frac \beta {
\sqrt{R_{\max }(t)}} \qquad \mbox{for}\;t \in [0,T_{\max }). 
\end{equation}
Then, by applying a result of Hamilton 
(see Theorems 16.4 and 16.5 in \cite{Ha3}), we know that there 
exists a sequence of dilations of the solution which 
converges to one of the singularity models of Type I, II or III.
We will analyse this limit in the next section.

We conclude this section with the following lemma 
which will be very useful in our analysis of the
Type I and Type II limits.

\vskip 3mm {\bf \underline{Lemma 3.1}} \quad 
%Let $\widetilde{g}_{\alpha \overline{\beta }}(\cdot ,t)$, 
%$(-\infty <t\leq 0)$ be a family of K\"ahler metric of 
%a complete two--dimensional manifold $
%\widetilde{M}$. 
Suppose 
$(\widetilde{M},\widetilde{g}_{\alpha \overline{\beta }}(\cdot ,t))$ 
is a complete ancient solution to the Ricci flow on a noncompact 
K\"ahler surface with nonnegative and bounded holomorphic 
bisectional curvature for all time. Then the curvature 
operator of the metric 
$\widetilde{g}_{\alpha \overline{\beta }}(\cdot ,t)$ 
is nonnegative definite everywhere on 
$\widetilde{M} \times (-\infty ,0]$.

\vskip 3mm {\bf \underline{Proof}} \quad 
Choose a local orthonormal coframe 
$\{\omega _1,\omega _2,\omega _3,\omega _4\}$ 
on an open set 
$U\subset \widetilde{M}$ so that 
$\omega _1+\sqrt{-1}\omega _2$ and $\omega_3+\sqrt{-1}\omega _4$ 
are $(1,0)$ forms over $U$. Then the self--dual forms 
$$
\varphi _1=\omega _1\land \omega _2+\omega _3\land \omega _4, 
\,\,\, \varphi_2=\omega _2\land \omega _3+\omega _1\land \omega _4, 
\,\,\, \varphi _3=\omega_3\land \omega _1+\omega _2\land \omega _4
$$ 
and the anti--self--dual forms 
$$
\psi_1=\omega _1\land \omega _2-\omega _3\land \omega _4, 
\,\,\, \psi _2=\omega_2\land \omega _3-\omega _1\land \omega _4, 
\,\,\, \psi _3=\omega _3\land \omega_1-\omega _2\land \omega _4
$$ 
form a basis of the space of $2$ forms over $U$.
In particular, $\varphi _1,\psi _1,\psi _2$ and $\psi _3$ 
give a basis for the space of $(1,1)$ forms over $U$.

On a K\"ahler surface, it is well known that its curvature 
operator has image in the holonomy algebra $u(2) \,(\subset so(4))$ 
spanned by $(1,1)$ forms. Thus, the curvature operator ${\bf M}$ 
in the basis $\{\varphi_1,\varphi _2,\varphi _3,\psi _1,\psi _2,\psi _3\}$ 
has the following form,
$$
{\bf M}=\left( 
\begin{array}{cc}
\begin{array}{ccc}
a & 0 & 0 \\ 
0 & 0 & 0 \\ 
0 & 0 & 0 
\end{array}
& 
\begin{array}{ccc}
b_1 & b_2 & b_3 \\ 
0 & 0 & 0 \\ 
0 & 0 & 0 
\end{array}
\\ 
\begin{array}{ccc}
b_1 & 0 & 0 \\ 
b_2 & 0 & 0 \\ 
b_3 & 0 & 0 
\end{array}
& A 
\end{array}
\right) \ 
$$
where $A$ is a $3\times 3$ symmetric matrix.

Let $V$ be a real tangent vector of the K\"ahler surface 
$\widetilde{M}$. Denote by $J$ the complex structure of the 
K\"ahler surface $\widetilde{M}$. It is clear that the complex 
$2$--plane $V\wedge JV$ is dual to $(1,1)$ form 
$u\varphi _1+v_1\psi _1+v_2\psi _2+v_3\psi _3$ satisfying the
decomposability condition $u^2=v_1^2+v_2^2+v_3^2$. Then after
normalizing $u$ to 1 by scaling, we see that the holomorphic 
bisectional curvature is nonnegative if and only if 
\begin{equation}
\label{3.4}
a+b\cdot v+b\cdot w+^{\ t}\negthinspace {}vAw\geq 0\ , 
\end{equation}
for any unit vectors $v=(v_1,v_2,v_3)$ and $w=(w_1,w_2,w_3)$ 
in ${\bf R^3}$, where $b$ is the vector $(b_1,b_2,b_3)$ 
in $\bf{M}$.

Denote by $a_1\leq a_2\leq a_3$ the eigenvalues of $A$. 
Recall that $\mbox{tr}\,A=a $ by the Bianchi identity, 
so if we choose $v$ to be the eigenvector of $A $ with
eigenvalue $a_3$ and choose $w=-v$, (\ref{3.4}) gives 
\begin{equation}
\label{3.5}a_1+a_2\geq 0\ .
\end{equation}
In particular, we have $a_2\geq 0$.

To proceed further, we need to adapt the maximum principle for parabolic
equations on compact manifold in Hamilton \cite{Ha1} to $\widetilde{M}$.
%In \cite{Ha1} Hamilton obtain a useful maximum principle for parabolic
%equations on compact manifolds.
Let
$$
\left( a_i\right) _{\min }(t)=\inf\limits_{x\in \widetilde{M}}a_i(x,t)\
,\qquad i=1,2,3\ 
$$
and
$$
K=\sup \limits_{(x,t)\in \widetilde{M}\times (-\infty ,0]}\left|
Rm(x,t)\right|.
$$
By assumption, the ancient solution
$\widetilde{g}_{\alpha \overline{\beta }}(\cdot ,t)$
has bounded holomorphic bisectional curvature, hence $K$ is finite.
Thus, by the derivative estimate of Shi \cite{Sh1} 
(see also Theorem 7.1 in \cite{Ha3}), the higher order 
derivatives of the curvature are also uniformly bounded. 
In particular, we can use the maximum
principle of Cheng--Yau (see Proposition 1.6 in \cite{CY}) and then, as
observed in \cite{Ha3}, this implies that the maximum principle of 
Hamilton in \cite{Ha1} actually works for the evolution equations 
of the curvature of $\widetilde{g}_{\alpha \overline{\beta }}(\cdot ,t)$ 
on the complete noncompact manifold $\widetilde{M}.$ Thus, from \cite{Ha1},
we obtain
\begin{eqnarray}
\frac{d\left( a_1\right) _{\min }}{dt}&\geq&\left( \left( a_1\right)_
{\min }\right) ^2+2\left( a_2\right) _{\min }\cdot 
\left( a_3\right) _{\min } \nonumber \\
& \geq & 3\left( \left( a_1\right) _{\min }\right) ^2 \nonumber
\end{eqnarray}
by (\ref{3.5}). Then, for fixed $t_0\in (-\infty ,0)$ and 
$t>t_0,$
\begin{eqnarray}
\left( a_1\right) _{\min }(t)&\geq&\frac 1{\left( a_1\right)_
{\min }^{-1}\left( t_0\right) -3\left( t-t_0\right) } \nonumber \\
& \geq & \frac 1{-K^{-1}-3\left( t-t_0\right) }\ . \nonumber
\end{eqnarray}
Letting $t_0\rightarrow -\infty $, we get 
\begin{equation}
\label{3.6}
a_1\geq 0\ ,\qquad \mbox{for all} \;\;\; (x,t) \in 
\widetilde{M} \times (-\infty ,0]\,
\end{equation}
i.e. $A \geq 0$.

Finally, to prove the nonnegativity of the curvature operator 
$\bf{M}$, we recall its corresponding ODE from \cite{Ha1},
%the corresronding ODE for the curvature 
%operator ${\bf M}$ is given by
$$
\frac{d{\bf M}}{dt}={\bf M}^2+\left( 
\begin{array}{cc}
0 & 0 \\ 
0 & A^{\#} 
\end{array}
\right) \ , 
$$
where $A^{\#} \geq 0$ is the adjoint matrix of $A.$

Let $m_1$ be the smallest eigenvalue of the curvature operator ${\bf M}$. 
By using the maximum principle of Hamilton ( \cite{Ha1} or \cite
{Ha3} ) again, we have
$$
\frac{d\left( m_1\right) _{\min }}{dt}\geq \left( m_1\right) _{\min }^2 
$$
where $\left( m_1\right) _{\min }(t)=\inf \limits_{x\in \widetilde{M}%
}m_1(x,t)$. Therefore, by the same reasoning in the derivation of 
(\ref{3.6}), we have
\begin{equation}
\label{3.7}
m_1\geq 0\ ,\qquad \mbox{for all} \;\; (x,t)\in \widetilde{M}\times
(-\infty ,0]\ . 
\end{equation}
So ${\bf{M}} \geq 0$ and the proof of the lemma is completed.
\hfill$\Box $

\section*{\S4. Time decay estimate on curvature}

\setcounter{section}{4} \setcounter{equation}{0}

\qquad
Let $(M,g_{\alpha \overline{\beta }}(x))$ be a complete
noncompact K\"ahler surface satisfying all the assumptions in
the Main Theorem and
$(M,g_{\alpha \overline{\beta }}(\cdot ,t))$, 
$t\in [0,T_{\max })$ be the maximal solution of the Ricci 
flow (\ref{2.1}) with $g_{\alpha \overline{\beta }}(\cdot )$ 
as the initial metric. Clearly, the maximal solution is of either
one of the following types.

\begin{description}
\item[{Type I:}]  \qquad $T_{\max }<+\infty $ and 
$\sup \left( T_{\max }-t\right) R_{\max }(t)<+\infty $;

\item[{Type II(a):}]  $T_{\max }<+\infty $ and $\sup 
\left( T_{\max}-t\right) R_{\max }(t)=+\infty $;

\item[{Type II(b):}]  $T_{\max }=+\infty $ and 
$\sup tR_{\max }(t)=+\infty $;

\item[{Type III:}]  \quad $T_{\max }=+\infty $ and 
$\sup tR_{\max }(t)<+\infty $.
\end{description}

In Section 3, we have proved that the maximal solution satisfies 
the following injectivity radius estimate
$$
\mbox{inj}(M,g_{\alpha \overline{\beta }}(\cdot ,t))\geq 
\frac \beta {\sqrt{R_{\max}(t)}}\ \qquad \mbox{on} \quad [0,T_{\max }), 
$$
for some $\beta >0$. By applying a result of Hamilton 
(Theorems 16.4 and 16.5 in \cite{Ha3}), we know that there exists 
a sequence of dilations of the solution converging to a singularity
model of the corresponding type. Note that since the maximal 
solution is complete and noncompact, the limit must also be 
complete and noncompact. The following is the main result of this 
section which says that this limit must be of Type III or equivalently, 
the maximal solution must be of Type III.
%The main purpose of this section is to prove the following time decay
%estimate for the solution of (\ref{2.1}).

\vskip 3mm{\bf \underline{Theorem 4.1}} 
\quad 
Let $(M,g_{\alpha \overline{\beta }}(x))$ be a complete
noncompact K\"ahler surface as above. 
%with bounded and positive holomorphic
%bisectional curvature. Suppose $(M,g_{\alpha \overline{\beta }}(x))$
%satisfies the condition $(i)$ and $(ii)^{\prime \prime }$. 
Then the Ricci flow (\ref{2.1}) with $g_{\alpha \overline{\beta }}(x)$ as 
the initial metric has a solution $g_{\alpha \overline{\beta }}(x,t)$ 
for all $t\in [0,+\infty) $ and $x\in M$. Moreover, the scalar 
curvature $R(x,t)$ of the solution satisfies 
\begin{equation}
\label{4.1}
0\leq R(x,t)\leq \frac C{1+t}\ \qquad \mbox{on} 
\quad M\times [0,+\infty)\ , 
\end{equation}
for some positive constant $C$.

\vskip 3mm{\bf \underline{Proof.}} \quad 
We prove by contradiction. Thus, suppose the maximal solution 
is of Type I or Type II and let 
$(\widetilde{M},\widetilde{g}_{\alpha \overline{\beta }}(x,t))$ 
be the limit of a sequence of dilations of the maximal solution which
%$( \widetilde{M},\widetilde{g}_{\alpha \overline{\beta }}(x,t))$ 
is then a singularity model of Type I or Type II respectively.
After a study of its properties, we can blow down the 
singularity model and apply a dimension 
reduction argument to obtain the desired contradiction. 

Now, recall that the maximal solution satisfies the maximal 
volume growth condition (i) by Proposition 2.1. Since condition 
(i) is also invariant under rescaling, we see that the singularity 
model $(\widetilde{M},\widetilde{g}_{\alpha \overline{\beta }}(x,t))$ 
also satisfies the maximal volume growth condition, i.e. 
\begin{equation}
\label{4.2}
\mbox{Vol}_t\left( \widetilde{B}_t(x,r)\right) \geq C_1r^4\ 
\qquad \mbox{for all} \quad 0\leq r< +\infty \quad \mbox{and} 
\quad x\in \widetilde{M}\ , 
\end{equation}
where $\mbox{Vol}_t\left( \widetilde{B}_t(x,r)\right) $ denotes 
the volume of the geodesic ball $\widetilde{B}_t(x,r)$ of radius $r$ with 
center at $x$ with respect to the metric 
$\widetilde{g}_{\alpha \overline{\beta }}(\cdot ,t).$

It is clear that the limit 
$\widetilde{g}_{\alpha \overline{\beta }}(\cdot ,t)$
has nonnegative holomorphic bisectional curvature. Thus, 
from Lemma 3.1, the curvature operator of the metric 
$\widetilde{g}_{\alpha \overline{\beta }}(\cdot ,t)$ 
is nonnegative definite everywhere.

Denote by $\widetilde{R}(\cdot ,t)$ the scalar curvature of 
$\widetilde{g}_{\alpha \overline{\beta }}(\cdot ,t)$ 
and $\widetilde{d}_t(x,x_0)$ the geodesic distance between two points 
$x,x_0\in \widetilde{M}$ with respect to the metric 
$\widetilde{g}_{\alpha \overline{\beta }}(\cdot,t)$. 
We claim that at time $t=0,$ we have
\begin{equation}
\label{4.3}
\limsup \limits_{\widetilde{d}_0(x,x_0)\rightarrow +\infty } 
\widetilde{R}(x,0)\widetilde{d}_0^2(x,x_0)=+\infty \  
\end{equation}
for any fixed $x_0\in \widetilde{M}.$

Suppose not, that is the curvature of the metric 
$\widetilde{g}_{\alpha \overline{\beta }}(\cdot ,0)$ 
has quadratic decay. Now, by applying a result of Shi 
(see Theorem 8.2 in \cite{Sh4}), we know that the solution 
$\widetilde{g}_{\alpha \overline{\beta }}(\cdot ,t)$ of the Ricci flow 
exists for all $t\in (-\infty ,+\infty )$ and satisfies 
\begin{equation}
\label{4.4}
\lim \limits_{t\rightarrow +\infty }\sup \left\{ \left. 
\widetilde{R}(x,t)\right| \ x\in \widetilde{M}\right\} =0\ . 
\end{equation}
On the other hand, by the Harnack inequality of Cao \cite{Cao}, 
we have 
\begin{equation}
\label{4.5} 
\frac{\partial \widetilde{R}}{\partial t}\geq 0\ \qquad \mbox{on} 
\quad \widetilde{M}\times (-\infty ,+\infty )\ . 
\end{equation}
Thus, combining (\ref{4.4}) and (\ref{4.5}), we deduce that
$$
\widetilde{R}\equiv 0\ ,\qquad \mbox{for all} \quad (x,t)\in 
\widetilde{M} \times (-\infty ,+\infty )\ 
$$
and hence $\widetilde{g}_{\alpha \overline{\beta }}(\cdot ,t)$ 
is flat for all $t\in (-\infty ,+\infty )$. But, by definition, 
a singularity model cannot be flat. This proves our claim (\ref{4.3}).

With the estimate (\ref{4.3}), we can then apply a lemma of Hamilton 
(Lemma 22.2 in \cite{Ha3}) to find a sequence of positive numbers 
$\delta _j$, $j=1,2,\cdots$, with $\delta_j\rightarrow 0$ such that

\begin{enumerate}
\item[{(a)}]  $\widetilde{R}(x,0)\leq (1+\delta _j)\widetilde{R}(x_j,0)$ 
for all $x$ in the ball $\widetilde{B}_0(x_j,r_j)$ of radius $r_j$ 
centered at $x_j$ with respect to the metric 
$\widetilde{g}_{\alpha \overline{\beta }}(\cdot,0)$;

\item[{(b)}]  $r_j^2\widetilde{R}(x_j,0)\rightarrow +\infty $;

\item[{(c)}]  if $s_j=\widetilde{d}_0(x_j,x_0)$, then $\lambda
_j=s_j/r_j\rightarrow +\infty $;

\item[{(d)}]  the balls $\widetilde{B}_0(x_j,r_j)$ are disjoint.
\end{enumerate}

Denote the minimum of the holomorphic sectional curvature of the metric $
\widetilde{g}_{\alpha \overline{\beta }}(\cdot ,0)$ at $x_j$ by $h_j$. We
claim that the following holds 
\begin{equation}
\label{4.6}
\varepsilon _j=\frac{h_j}{\widetilde{R}(x_j,0)}\rightarrow 0\
\qquad \mbox{as} \quad j\rightarrow +\infty \ . 
\end{equation}

Suppose not, there exists a subsequence $j_k\rightarrow +\infty $ 
and some positive number $\varepsilon >0$ such that 
\begin{equation}
\label{4.7}
\varepsilon _{j_k}=\frac{h_{j_k}}{\widetilde{R}(x_{j_k},0)}\geq
\varepsilon \ \qquad \mbox{for all} \quad k=1,2,\cdots \ . 
\end{equation}

Since the solution $\widetilde{g}_{\alpha \overline{\beta }}(\cdot ,t)$
is ancient, it follows from the Harnack inequality of Cao \cite{Cao} that
the scalar curvature $\widetilde{R}(x,t)$ is pointwisely nondecreasing in
time. Then, by using the local derivative estimate of Shi \cite{Sh1} 
(or see Theorem 13.1 in \cite{Ha3}) and (a), (b), we have
\begin{eqnarray}
\label{4.8}
\sup \limits_{x\in \widetilde{B}_0(x_{j_k},r_{j_k})}\left| 
\nabla \widetilde{R}m(x,0)\right| ^2&\leq&C_5\widetilde{R}^2(x_j,0)
\left( \frac 1{r_{j_k}^2}+\widetilde{R}(x_j,0)\right) \nonumber \\
&\leq&2C_5\widetilde{R}^3(x_j,0)\ ,
\end{eqnarray}
where $\widetilde{R}m$ is the curvature tensor of 
$\widetilde{g}_{\alpha \overline{\beta }}$ and $C_5$ is a 
positive constant depending only on the dimension.

For any $x\in \widetilde{B}_0(x_{j_k},r_{j_k})$, we obtain 
from (\ref{4.7}) and (\ref{4.8}) that the minimum of the 
holomorphic sectional curvature $h_{\min}(x)$ at $x$, satisfies
\begin{eqnarray}
\label{4.9}
h_{\min }(x)&\geq&h_{j_k}-\sqrt{2C_5}\widetilde{R}^{3/2}(x_{j_k},0)\widetilde{d}_0(x,x_{j_k})\nonumber\\
&\geq&\widetilde{R}(x_{j_k},0)\left( \varepsilon -\sqrt{2C_5}\cdot \sqrt{\widetilde{R}(x_{j_k},0)}\cdot \widetilde{d}_0(x,x_{j_k})\right)\nonumber\\
&\geq&\frac \varepsilon 2\widetilde{R}(x_{j_k},0)\ 
\end{eqnarray}
if
$$
\widetilde{d}_0(x,x_{j_k})\leq \frac{\varepsilon}{2\sqrt{2C_5}\cdot \sqrt{ 
\widetilde{R}(x_{j_k},0)}}\ . 
$$
Thus, from (a) and (\ref{4.9}), there exists $k_0>0$ such that
for any $k\geq k_0$ and 
$$
x \in \widetilde{B}_0(x_{j_k},\frac{\varepsilon}{2\sqrt{2C_5}
\cdot \sqrt{\widetilde{R}(x_{j_k},0)}}),
$$
we have
\begin{equation}
\label{4.10}
\frac \varepsilon 2\widetilde{R}(x_{j_k},0) \leq 
\mbox{holomorphic sectional curvature at $x$} \leq 
2\widetilde{R}(x_{j_k},0). 
\end{equation}
%in the geodesic ball $\widetilde{B}_0(x_{j_k},\frac \varepsilon {2\sqrt{2C_5}%
%\cdot \sqrt{\widetilde{R}(x_{j_k},0)}})$.

We have proved that the metric 
$\widetilde{g}_{\alpha \overline{\beta }}(\cdot ,0)$ 
has nonnegative definite curvature operator. In particular, the
sectional curvature is nonnegative. Then, by the generalized 
Cohn--Vossen inequality in real dimension 4 \cite{GW1}, we have 
\begin{equation}
\label{4.11}
\int_{\widetilde{M}}\Theta \leq \chi \left( \widetilde{M}\right)
<+\infty 
\end{equation}
where $\Theta $ is the Gauss--Bonnet--Chern integrand for 
the metric $\widetilde{g}_{\alpha \overline{\beta }}(\cdot ,0)$ 
and $\chi \left( \widetilde{M}\right)$ is the Euler number of the 
manifold $\widetilde{M}$ which has
finite topology type by the soul theorem of Cheeger--Gromoll.

On the other hand, from the proof of Theorem 1.3 of Bishop--Goldberg 
\cite{BG} (see Page 523 of \cite{BG}), the inequality (\ref{4.10}) 
implies that 
\begin{equation}
\label{4.12}
\Theta (x)\geq C(\varepsilon)\widetilde{R}^2(x_{j_k},0)\ 
\qquad \mbox{for all} \;\;\; x\in \widetilde{B}_0(x_{j_k},
\frac{\varepsilon}{2\sqrt{2C_5}\cdot \sqrt{ \widetilde{R}(x_{j_k},0)}}), 
\end{equation}
%for all $x\in \widetilde{B}_0(x_{j_k},\frac \varepsilon{2\sqrt{2C_5}\cdot \sqrt{ 
%\widetilde{R}(x_{j_k},0)}})$, 
where $C(\varepsilon)$ is some positive constant depending only on 
$\varepsilon$. Now, by combining (\ref{4.2}), (b), (d), (\ref{4.11})
and (\ref{4.12}), we get
\begin{eqnarray*}
+\infty >\chi \left( \widetilde{M}\right)
& \geq & \sum\limits_{k=k_0}^\infty \int_{\widetilde{B}_0(x_{j_k},
\frac{\displaystyle \varepsilon}{2\sqrt{2C_5}\cdot 
\sqrt{ \widetilde{R}(x_{j_k},0)}})}\Theta \\
& \geq & C(\varepsilon)\sum\limits_{k=k_0}^\infty \widetilde{R}^2
(x_{j_k},0)\cdot C_1\left( \frac\varepsilon {2\sqrt{2C_5}\cdot \sqrt{\widetilde{R}(x_{j_k},0)}}\right) ^4 \\
& = & C(\varepsilon )\sum\limits_{k=k_0}^\infty 
\frac{C_1\varepsilon ^4}{64C_5^2}\\
& = & +\infty
\end{eqnarray*}
which is a contradiction. Hence our claim (\ref{4.6}) is proved.

Now, we are going to blow down the singularity model $(\widetilde{M}, 
\widetilde{g}_{\alpha \overline{\beta }}(x,t))$. For the above chosen $x_j$, 
$r_j$ and $\delta _j$, let $x_j$ be the new origin $O$, dilate the space by
a factor $\lambda _j$ so that $\widetilde{R}(x_j,0)$ become $1$ at the
origin at $t=0$, and dilate in time by $\lambda _j^2$ so that it is still a
solution to the Ricci flow. The balls $\widetilde{B}_0(x_j,r_j)$ are dilated
to the balls centered at the origin of radii $\widetilde{r}_j=r_j^2 
\widetilde{R}(x_j,0)\rightarrow +\infty $ ( by (b) ). Since the scalar
curvature of $\widetilde{g}_{\alpha \overline{\beta }}(x,t)$ is pointwise
nondecreasing in time by the Harnack inequality, the curvature bounds on $
\widetilde{B}_0(x_j,r_j)$ also give bounds for previous times in these
balls. And the maximal volume growth estimate (\ref{4.2}) and the local
injectivity radius estimate of Cheeger, Gromov and Taylor \cite{CGT} imply
that
$$
\mbox{inj}_{\widetilde{M}}\left( x_j,\widetilde{g}_{\alpha 
\overline{\beta }}(\cdot,0)\right) 
\geq \frac \beta {\sqrt{\widetilde{R}(x_j,0)}}\ ,
$$
for some positive constant $\beta $ independent of $j.$

So we have everything to take a limit for the dilated solutions. By applying
the compactness theorem in \cite{Ha2} and combining (\ref{4.2}), (\ref{4.6}%
), (a) and (b), we obtain a complete noncompact solution, still denoted by $%
( \widetilde{M},\widetilde{g}_{\alpha \overline{\beta }}(x,t))$, for $t\in
(-\infty ,0]$ such that

\begin{enumerate}
\item[{(e)}]  the curvature operator is still nonnegative;

\item[{(f)}]  $\widetilde{R}(x,t)\leq 1$, for all $x\in \widetilde{M}$, $
t\in (-\infty ,0]$, and $\widetilde{R}(0,0)=1$;

\item[{(g)}]  $\mbox{Vol}_t\left( \widetilde{B}_t(x,r)\right) \geq C_1r^4$ 
for all $x\in \widetilde{M}$, $0\leq r\leq +\infty $;

\item[{(h)}]  there exists a complex $2$--plane $V\wedge JV$ at the origin $%
O $ so that at $t=0$, the corresponding holomorphic sectional curvature
vanishes.
\end{enumerate}

If we consider the universal covering of $\widetilde{M}$, the induced metric
of $\widetilde{g}_{\alpha \overline{\beta }}(\cdot ,t)$ on the universal
covering is clearly still a solution to the Ricci flow and satisfies all 
of above (e), (f), (g), (h). Thus, without loss of generality, we may assume 
that $\widetilde{M}$ is simply connected.

Next, by using the strong maximum principle on the evolution equation of the
curvature operator of $\widetilde{g}_{\alpha \overline{\beta }}(\cdot ,t)$
as in \cite{Ha1} (see Theorem 8.3 of \cite{Ha1}), we know that there
exists a constant $K>0$ such that on the time interval $-\infty <t<-K$, the
image of the curvature operator of $(\widetilde{M},\widetilde{g}_{\alpha 
\overline{\beta }}(\cdot ,t))$ is a fixed Lie subalgebra of $so(4)$ of
constant rank on $\widetilde{M}$. Because $\widetilde{M}$ is K\"ahler, the
possibilities are limited to $u(2)$, $so(2)\times so(2)$ or $so(2).$

In the case $u(2)$, the sectional curvature is strictly positive. Thus, this
case is ruled out by (h). In the cases $so(2)\times so(2)$ or $so(2)$,
according to \cite{Ha1}, the simply connected manifold $\widetilde{M}$
splits as a product $\widetilde{M}=\Sigma _1\times \Sigma _2$, where $\Sigma
_1$ and $\Sigma _2$ are two Riemann surfaces with nonnegative curvature (by
(e)), and at least one of them, say $\Sigma _1$, has positive curvature (by (f)).

Denote by $\widetilde{g}_{\alpha \overline{\beta }}^{(1)}(\cdot ,t)$ the
corresponding metric on $\Sigma _1$. Clearly, it follows from (g) and
standard volume comparison that for any $x\in \Sigma _1$, $t\in (-\infty
,-K) $, we have 
\begin{equation}
\label{4.13}
\mbox{Vol}B_{\Sigma_1}(x,r) \geq C_6 r^2 \qquad \mbox{for} \;\; 
0 \leq r  < +\infty
\end{equation}
%\left\{ 
%\begin{array}{c}
%the\ volume\ of\ the\ geodesic\ ball\ of\ \Sigma _1 \\ 
%with\ center\ at\ x\ and\ radius\ r \\ 
%with\ respect\ the\ metric\ \widetilde{g}_{\alpha \overline{\beta }%
%}^{(1)}(\cdot ,t) 
%\end{array}
%\right\} \geq C_6r^2\ ,\qquad for\quad 0\leq r<+\infty. 
%\end{equation}
where both the geodesic ball $B_{\Sigma_1}(x,r)$ 
and the volume are taken with respect to the metric 
$\widetilde{g}_{\alpha \overline{\beta }}^{(1)}(\cdot ,t)$ on $\Sigma_1$,
$C_6$ is a positive constant depending only on $C_1.$
Also as the curvature of 
$\widetilde{g}_{\alpha \overline{\beta }}^{(1)}(x,t)$ 
is positive, it follows from Cohn--Vossen inequality that 
\begin{equation}
\label{4.14}
\int_{\Sigma _1}\widetilde{R}^{(1)}(x,t)d\sigma _t\leq 8\pi \ , 
\end{equation}
where $\widetilde{R}^{(1)}(x,t)$ is the scalar curvature of $(\Sigma _1, 
\widetilde{g}_{\alpha \overline{\beta }}^{(1)}(x,t))$ and $d\sigma _t$ is
the volume element of the metric $\widetilde{g}_{\alpha \overline{\beta }%
}^{(1)}(x,t).$

Now, the metric 
$\widetilde{g}_{\alpha \overline{\beta }}^{(1)}(x,t)$ 
is a solution to the Ricci flow on the Riemann surface $\Sigma _1$ 
over the ancient time interval $(-\infty ,-K)$. 
Thus, (\ref{4.13}) and (\ref{4.14}) imply that for each $t\in (-\infty ,-K)$,
the curvature of $\widetilde{g}_{\alpha \overline{\beta }}^{(1)}(x,t)$ has
quadratic decay in the average sense of Shi \cite{Sh4} and then the apriori
estimate of Shi (see Theorem 8.2 in \cite{Sh4}) implies that the solution 
$\widetilde{g}_{\alpha \overline{\beta }}^{(1)}(x,t)$ exists for all 
$t\in (-\infty ,+\infty )$ and satisfies
\begin{equation}
\label{4.15}
\lim \limits_{t\rightarrow +\infty }\sup \left\{ \left. 
\widetilde{R}^{(1)}(x,t)\right| \ x\in \Sigma _1\right\} =0.
\end{equation}
Again, by the Harnack inequality of Cao \cite{Cao}, 
we know that $\widetilde{R}^{(1)}(x,t)$ is pointwisely 
nondecreasing in time. Therefore, we conclude that
$$
\widetilde{R}^{(1)}(x,t)\equiv 0\qquad \mbox{on}\quad 
\Sigma _1\times (-\infty,+\infty )\ . 
$$
This contradicts with the fact that $(\Sigma _1,\widetilde{g}_{\alpha 
\overline{\beta }}(\cdot ,t))$ has positive curvature for $t<-K.$ Hence we
have seeked the desired contradiction and have completed the 
proof of Theorem 4.1.
\hfill $\Box $

\section*{\S5. Topology and Steinness}

\setcounter{section}{5} \setcounter{equation}{0}

\qquad 
In this section, we use the estimates obtained in the previous
sections to study the topology and the complex structure of the
K\"ahler surface in our Main Theorem. Our result is

\vskip 3mm{\bf \underline{Theorem 5.1}} \quad 
Suppose $(M,g_{\alpha \overline{\beta }})$ is a complete 
noncompact K\"ahler surface satisfying the assumptions in the
Main Theorem. 
%with bounded and positive holomorphic bisectional curvature and satisfies
%the condition $(i)$ and $(ii)^{\prime \prime }$. 
Then $M$ is homeomorphic to $\R^4$ and is a Stein manifold.

\vskip 3mm 
The proof of this theorem is exactly the same as
in \cite{CZ2}. For the convenience of the readers, we give a 
sketch of the arguments and refer to the cited reference for details.

\vskip 3mm{\bf \underline{Sketch of proof.}} \qquad 
We evolve the
metric $g_{\alpha \overline{\beta }}(x)$ by the Ricci flow (\ref{2.1}). From
Theorem 4.1, the solution $g_{\alpha \overline{\beta }}(x,t)$ exists for all $
t\in [0,+\infty )$ and satisfies 
\begin{equation}
\label{5.1}
R(x,t)\leq \frac C{1+t}\ \qquad \mbox{on} \quad M\times [0,+\infty )\ 
\end{equation}
for some positive constant $C$. Also, Proposition 2.1 tells us 
that the volume growth condition (i) is preserved under 
the Ricci flow. By using the
local injectivity radius estimate of Cheeger--Gromov--Taylor, 
this implies that 
\begin{equation}
\label{5.2}
\mbox{inj}\left( M,g_{\alpha \overline{\beta }}(\cdot ,t)\right) \geq
C_7(1+t)^{\frac 12}\ \qquad \mbox{for}\quad t\in [0,+\infty )\  
\end{equation}
with some positive constant $C_7$.

Since the Ricci curvature of $g_{\alpha \overline{\beta }}(x,t)$ 
is positive for all $x\in M$ and $t\geq 0$, the Ricci flow equation 
(\ref{2.1}) implies that the ball $B_t(x_0,\frac{C_7}2(1+t)^{\frac 12})$
of radius $\frac{C_7}2(1+t)^{\frac 12}$ with respect to the metric 
$g_{\alpha \overline{\beta }}(\cdot ,t)$ contains the ball 
$B_0(x_0,\frac{C_7}2(1+t)^{\frac 12})$ of the same radius with 
respect to the initial metric $g_{\alpha \overline{\beta }}(\cdot ,0)$. 
Combining this with (\ref{5.2}), we deduce that
$$
\pi _p(M,x_0)=0\ \qquad \mbox{for any} \quad p\geq 1\ 
$$
and
$$
\pi _q(M,\infty )=0\ \qquad \mbox{for} \quad 1\leq q\leq 2\ , 
$$
where $\pi _q(M,\infty )$ is the $q$th homotopy group of $M$ at infinity.

Thus, by the resolution of the generalized Poincar\'e conjecture 
on four manifolds by Freedman \cite{Fre}, we know that 
$M$ is homeomorphic to $\R^4$.

Next, the injectivity radius estimate (\ref{5.2}) also tells us that, 
for $t$ large enough, the exponential maps provide diffeomorphisms
between big geodesic balls $B_t(x_0,\frac{C_7}2(1+t)^{\frac 12})$ 
of $M$ with big Euclidean balls on $\C^2$. The curvature estimate 
(\ref{5.1}) together with its derivative estimates imply that the 
difference of the complex structure of $M$ and the standard complex 
struture of $\C^2$ in those geodesic balls can be made arbitrarily small 
by taking $t$ large enough. Then, we can use the $L^2$
estimates of H\"ormander \cite{Ho} to modify the exponential map 
to a biholomorphism between a domain $\Omega (t)$
containing $B_t(x_0,\frac{C_7}4(1+t)^{\frac 12})$ and a Euclidean ball of
radius $\frac{C_7}3(1+t)^{\frac 12}$ in $\C^2$. Since the solution 
$g_{\alpha \overline{\beta }}(\cdot ,t)$ of the Ricci flow is shrinking, 
the domain $\Omega (t)$ contains the geodesic ball 
$B_0(x_0,\frac{C_7}4(1+t)^{\frac 12})$. Thus, we
can choose a sequence of $t_k\rightarrow +\infty $, such that
\begin{displaymath}
M=\cup _{k=0}^{+\infty }\Omega (t_k)\ ,\qquad \Omega (t_1)\subset \Omega
(t_2)\subset \cdots \subset \Omega (t_k)\subset \cdots \ . 
\end{displaymath}
Since for each $k$, $\Omega (t_k)$ is biholomorphic to 
the unit ball of $\C^2$, $(\Omega (t_k),\Omega (t_l))$ forms a Runge pair 
for any $k,\ l$.
Finally, we can appeal to a theorem of Markeo \cite{Mar} 
(see also Siu \cite{Si}) to conclude that $M$ is a Stein manifold.
\hfill $\Box $

\section*{\S6. Space decay estimate on curvature and the
Poincar\'e--Lelong equation}

\setcounter{section}{6} \setcounter{equation}{0}

\qquad 
Let $(M,g_{\alpha \overline{\beta }})$ be a complete noncompact
K\"ahler surface satisfying all the assumptions in the Main Theorem.
The main purpose of this section is to establish the existence of a strictly
plurisubharmonic function of logarithmic growth on $M$. To this end,
we first prove a curvature decay estimate at infinity of the metric 
$g_{\alpha \overline{\beta }}$.

\vskip 3mm{\bf \underline{Theorem 6.1}} \quad 
Let $(M,g_{\alpha \overline{\beta }})$ be a complete noncompact 
K\"ahler surface as above.
%with bounded and positive bisectional curvature. Suppose $(M,g_{\alpha 
%\overline{\beta }})$ satisfies the condition $(i)$ and $(ii)^{\prime \prime
%} $. 
Then there exists a constant $C>0$ such that for all $x\in M,\ r>0$, we
have 
\begin{equation}
\label{6.1}
\int_{B(x,r)}R(y)\,\frac 1{d^2(x,y)}\,dy \leq C\log (2+r)\ . 
\end{equation}

\vskip 3mm{\bf \underline{Proof.}} \quad 
Let $g_{\alpha \overline{\beta }}(x,t)$ be the solution of 
the Ricci flow (\ref{2.1}) with $g_{\alpha \overline{\beta }}(x)$ 
as the initial metric. From Theorem 4.1, we know that the 
solution exists for all times and satisfies 
\begin{equation}
\label{6.2}
R(x,t) \leq \frac{C_8}{1+t}\ \qquad \mbox{on} \quad 
M\times [0,+\infty ) 
\end{equation}
for some positive constant $C_8.$

Let 
$$
F(x,t)=\log \frac{\det \left( g_{\alpha \overline{\beta }}
(x,t)\right) }{\det \left( g_{\alpha \overline{\beta }}(x,0)\right)} 
$$
be the function introduced in the proof of Proposition 2.1. Since
$$
-\partial _\alpha \overline{\partial }_\beta \log \frac{\det \left(
g_{\gamma \overline{\delta }}(\cdot ,t)\right) }{\det \left( g_{\gamma 
\overline{\delta }}(\cdot ,0)\right) }=R_{\alpha \overline{\beta }}(\cdot
,t)-R_{\alpha \overline{\beta }}(\cdot ,0)\ , 
$$
after taking trace with the initial metric 
$g_{\alpha \overline{\beta }}(\cdot,0)$, we get 
\begin{equation}
\label{6.3}
R(\cdot ,0)=\bigtriangleup _0F(\cdot ,t)+g^{\alpha \overline{
\beta }}(\cdot ,0)R_{\alpha \overline{\beta }}(\cdot ,t) 
\end{equation}
where $\bigtriangleup _0$ is the Laplace operator of the metric 
$g_{\alpha \overline{\beta }}(\cdot ,0).$

Since $(M,g_{\alpha \overline{\beta }}(\cdot ,0))$ has positive Ricci
curvature and maximal volume growth, it is
well known (see \cite{ScY}) that the Green function $G_0(x,y)$ of
the initial metric $g_{\alpha \overline{\beta }}(\cdot ,0)$ exists on $M$
and satisfies the estimates 
\begin{equation}
\label{6.4}
\frac{C_9^{-1}}{d_0^2(x,y)}\leq G_0(x,y)\leq 
\frac{C_9}{d_0^2(x,y)} 
\end{equation}
and 
\begin{equation}
\label{6.5}
\left| \nabla _yG_0(x,y)\right| _0\leq \frac{C_9}{d_0^3(x,y)} 
\end{equation}
for some positive constant $C_9$ depending only on $C_1.$

For any fixed $\overline{x}_0\in M$ and any $\alpha >0$, we denote
$$
\Omega _\alpha =\left\{ \left. x\in M\right| \ G_0
(\overline{x}_0,x)\geq \alpha \ \right\} \ . 
$$
By (\ref{6.4}), it is not hard to see 
\begin{equation}
\label{6.6}
B_0\left( \overline{x}_0,\left( \frac{C_9^{-1}}\alpha \right)
^{\frac 12}\right) \subset \Omega _\alpha \subset B_0\left( \overline{x}
_0,\left( \frac{C_9}\alpha \right) ^{\frac 12}\right) \ . 
\end{equation}
Recall that $F$ evolves by 
$$
\frac{\partial F(x,t)}{\partial t}=-R(x,t)\ \qquad \mbox{on} \quad M\times
[0,+\infty )\ . 
$$
Combining with (\ref{6.2}), we obtain 
\begin{equation}
\label{6.7}0\geq F(x,t)\geq -C_{10}\log (1+t) \qquad \mbox{on} \quad M\times
[0,+\infty )\ . 
\end{equation}
Multiplying (\ref{6.3}) by $G_0(\overline{x}
_0,x)-\alpha $ and integrating over $\Omega _\alpha $, we have
\begin{eqnarray}
\label{6.8}
\int_{\Omega _\alpha }R(x,0)\left( G_0(\overline{x}_0,x)-\alpha \right) 
dx&=&\int_{\Omega _\alpha }\left( \bigtriangleup _0F(x,t)\right) 
\left( G_0(\overline{x}_0,x)-\alpha \right) dx \nonumber \\
&  & +\int_{\Omega _\alpha }g^{\alpha \overline{\beta }}(x,0)
R_{\alpha \overline{\beta }}(\cdot ,t)\left( G_0(\overline{x}_0,x)-
\alpha \right) dx  \nonumber \\
& = & -\int_{\partial \Omega _\alpha }F(x,t)\frac{\partial G_0
(\overline{x}_0,x)}{\partial \nu }d\sigma -F(\overline{x}_0,t) \nonumber\\
& & +\int_{\Omega _\alpha }g^{\alpha \overline{\beta }}(x,0)
R_{\alpha \overline{\beta }}(\cdot ,t)\left( G_0(\overline{x}_0,x)-\alpha 
\right) dx\nonumber\\
& \leq & C_{10}\left( 1+C_9^{\frac 52}\alpha ^{\frac 32}\mbox{Vol}_0
\left( \partial \Omega _\alpha \right) \right) \log (1+t) \nonumber \\
& & +\int_{\Omega _\alpha }g^{\alpha \overline{\beta }}(x,0)
R_{\alpha \overline{\beta }}(\cdot ,t)G_0(\overline{x}_0,x)dx\ ,\nonumber\\
\end{eqnarray}
by (\ref{6.4}) and (\ref{6.7}). Here, we have used $\nu $ to denote the outer unit
normal of $\partial \Omega _\alpha .$

From the coarea formula, we have
\begin{eqnarray}
\frac 1\alpha \int_\alpha ^{2\alpha }r^{\frac 32}Vol_0\left( 
\partial \Omega _r\right) dr&\leq&2^{\frac 32}\alpha ^{\frac 12}
\int_\alpha ^{2\alpha }\int_{\partial \Omega _r}\left| \nabla 
G_0(\overline{x}_0,x)\right| _0d\sigma \left| d\nu \right|  \nonumber \\
& \leq & 2^{\frac 32}C_9^{\frac 52}\alpha ^2\mbox{Vol}_0\left( 
\Omega _\alpha\right) \nonumber  \\
& \leq & 2^{\frac 32}C_9^{\frac 52}\alpha ^2\mbox{Vol}_0\left( 
B_0\left( \overline{x}_0,\left( \frac{C_9}\alpha \right) ^{\frac 12}\right) 
\right)\nonumber\\
& \leq & C_{11} \nonumber
\end{eqnarray}
for some positive constant $C_{11}$ by the standard volume
comparison. Substitute this into (\ref{6.8}) and
integrate (\ref{6.8}) from $\alpha $ to $2\alpha $, we get
\begin{eqnarray}
\label{6.9}
\int_{\Omega _{2\alpha}}R(x,0)\left( G_0(\overline{x}_0,x)-2\alpha 
\right) dx&\leq&C_{10}\left( 1+C_9^{\frac 52}C_{11}\right) \log (1+t)
\nonumber\\
&   &  +\int_{\Omega _\alpha }g^{\alpha \overline{\beta }}(x,0)
R_{\alpha \overline{\beta }}(x,t)G_0(\overline{x}_0,x)dx\ .\nonumber\\
\end{eqnarray}
It is easy to see that
$$
\int_{\Omega _{4\alpha }}R(x,0)G_0(\overline{x}_0,x)dx\leq 2\int_{\Omega
_{2\alpha }}R(x,0)\left( G_0(\overline{x}_0,x)-2\alpha \right) dx 
$$
and by the equation Ricci flow (\ref{2.1}), we also have
\begin{eqnarray}
&&\int_0^t\int_{\Omega _\alpha }g^{\alpha \overline{\beta }}(x,0)
R_{\alpha \overline{\beta }}(x,t)G_0(\overline{x}_0,x)dx dt 
\nonumber\\
& = & \int_{\Omega _\alpha }g^{\alpha \overline{\beta }}(x,0)
\left( g_{\alpha \overline{\beta }}(x,0)-g_{\alpha \overline{\beta }}
(x,t)\right) G_0(\overline{x}_0,x)dx \nonumber \\
& \leq & 2\int_{\Omega _\alpha }G_0(\overline{x}_0,x)dx\ .\nonumber
\end{eqnarray}
Thus by integrating (\ref{6.9}) in time from $0$ to $t$ and combining the
above two inequalities, we get for any $t>0,$
$$
\int_{\Omega _{4\alpha }}R(x,0)G_0(\overline{x}_0,x)dx\leq 2C_{10}\left(
1+C_9^{\frac 52}C_{11}\right) \log (1+t)+\frac 4t\int_{\Omega _\alpha }G_0( 
\overline{x}_0,x)dx\ . 
$$
Finally, substituting (\ref{6.4}) and (\ref{6.6}) into the above 
inequality, we see that there exists some positive constant $C_{12}$ 
such that for any $\overline{x}_0\in M$, $t>0$ and $r>0,$
\begin{equation}
\label{6.10}\int_{B_0(\overline{x}_0,r)}R(x,0)\frac 1{d^2(\overline{x}%
_0,x)}dx\leq C_{12}\left( \log (1+t)+\frac{r^2}t\right) \ . 
\end{equation}
Choose $t=r^2$, we get the desired estimate.
\hfill  $\Box $

\vskip 3mm 
Now we can use the estimate (\ref{6.1}) to solve the following 
Poincar\'e--Lelong equation on $M$
\begin{equation}
\label{6.11}
\sqrt{-1}\partial \overline{\partial }u= \mbox{Ric} \
\end{equation}
to get the strictly plurisubharmonic
function mentioned at the beginning of this section. 

As in \cite{MSY} or \cite{NST}, we first study the corresponding Poisson
equation on $M$
\begin{equation}
\label{6.12}\bigtriangleup u=R. 
\end{equation}
After we solve the Poisson equation (\ref{6.12}) with a solution of 
logarithmic growth, we will see that it is indeed a solution of
the Poincar\'e--Lelong equation with logarithmic growth.
 
To solve (\ref{6.12}), we first construct a family of approximate
solutions $u_r$ as follows.

For a fixed $x_0\in M$ and $r>0$, define $u_r(x)$ on $B(x_0,r)$ by
$$
u_r(x)=\int_{B(x_0,r)}\left( G(x_0,y)-G(x,y)\right) R(y)dy 
$$
where $G(x,y)$ is the Green function of the metric 
$g_{\alpha \overline{\beta }}$ on $M$. It is clear that
$$
u_r(x_0)=0 \qquad \mbox{and} \qquad \bigtriangleup u_r(x)=R(x) \quad \mbox{on} \quad B(x_0,r). 
$$
For $x\in B(x_0,\frac r2)$, we write
\begin{eqnarray}
u_r(x) & = & \left( \int_{B(x_0,r)\backslash B(x_0,2d(x,x_0))}+
\int_{B(x_0,2d(x,x_0))}\right) \left( G(x_0,y)-G(x,y)\right) 
R(y)dy \nonumber \\
& := &I_1+I_2\ . \nonumber
\end{eqnarray}
From (\ref{6.1}), we see that 
\begin{equation}
\label{6.13}
\left| I_2\right| \leq C_{13}\log \left( 2+d(x,x_0)\right) 
\qquad \mbox{on} \quad B(x_0,\frac r2) 
\end{equation}
for some positive constant $C_{13}$ independent of $x_0$, $x$ and $r.$

To estimate $I_1$, we get from (\ref{6.5}) that for $y\in B(x_0,r)\backslash
B(x_0,2d(x,x_0)),$%
\begin{eqnarray}
\left| G(x_0,y)-G(x,y)\right| &\leq&d(x,x_0)\cdot \sup 
\limits_{z\in B(x_0,d(x,x_0))}\left| \nabla _zG(z,y)\right|
\nonumber \\
& \leq & C_9d(x,x_0)\cdot \sup \limits_{z\in B(x_0,d(x,x_0))}
\frac 1{d^3(z,y)} \nonumber \\
&\leq&8C_9\frac{d(x,x_0)}{d^3(y,x_0)}\ .\nonumber
\end{eqnarray}
Thus by (\ref{6.1}), we have
\begin{eqnarray}
\label{6.14}
\left|I_1\right|&\leq&8C_9d(x,x_0)\int_{B(x_0,r)\backslash B(x_0,2d(x,x_0))}\frac{R(y)}{d^3(y,x_0)}dy \nonumber\\
&\leq&8C_9d(x,x_0)\sum\limits_{k=1}^\infty \frac 1{2^kd(x,x_0)}\cdot \int_{B(x_0,2^{k+1}d(x,x_0))\backslash B(x_0,2^kd(x,x_0))} 
\frac{R(y)}{d^2(y,x_0)}dy \nonumber \\
& \leq & 8C_9C\sum\limits_{k=1}^\infty \frac 1{2^k}\log \left( 2+2^{k+1}d(x,x_0)\right)\nonumber\\
& \leq & C_{14}\log \left( 2+d(x,x_0)\right) 
\end{eqnarray}
for some positive constant $C_{14}.$

Hence, by combining (\ref{6.13}) and (\ref{6.14}), we deduce 
\begin{equation}
\label{6.15}
\left| u_r(x)\right| \leq \left( C_{13}+C_{14}\right) \log
\left( 2+d(x,x_0)\right) 
\end{equation}
for any $r\geq 2d(x,x_0).$

On the other hand, by taking the derivative of $u_r(x)$, we get
\begin{eqnarray}
\label{6.16}
\left| \nabla u_r(x)\right|&\leq&C_9\int_M\frac{R(y)}
{d^3(x,y)}dy \nonumber \\
& \leq & C_9\left( \int\limits_{B(x,1)}\frac{R(y)}{d^3(x,y)}dy
+\sum\limits_{k=1}^\infty \frac 1{2^{k-1}}\int\limits_{B(x,2^k)
\backslash B(x,2^{k-1})}\frac{R(y)}{d^2(x,y)}dy\right) \nonumber\\
& \leq & C_9\left( C_{15}+\sum\limits_{k=1}^\infty 
\frac 1{2^{k-1}}C\log \left( 2+2^k\right) \right) \nonumber \\
& = & C_{16}\ .
\end{eqnarray}
Here, we have used (\ref{6.1}) and (\ref{6.5}), $C_{15}$ and 
$C_{16}$ are positive constants independent of $r$. 
Therefore, it follows from the Schauder theory of elliptic equations 
that there exists a sequence of $r_j\rightarrow +\infty $ 
such that $u_{r_j}(x)$ converges uniformly on
compact subset of $M$ to a smooth function $u$ satisfying 
\begin{equation}
\label{6.17}
\left\{ 
\begin{array}{ll}
\bigbreak u(x_0)=0\qquad \mbox{and} \qquad \bigtriangleup 
u=R &  \mbox{on} \;\; M\ , \\ 
\bigbreak |u(x)|\leq \left( C_{13}+C_{14}\right) \log \left(
2+d(x,x_0)\right)   & \mbox{for} \;\; x\in M\ , \\ 
|\nabla u(x)|\leq C_{16}\ & \mbox{for} \;\; x\in M\ . 
\end{array}
\right. 
\end{equation}
Thus, we have obtained a solution $u$ of logarithmic growth to the Poisson
equation (\ref{6.12}) on $M$. In the following we prove that $u$ is
actually a solution of the Poincar\'e--Lelong equation (\ref{6.11}).

Recall the Bochner identity, with $\bigtriangleup u = R$
\begin{eqnarray}
\label{6.18}
\frac 12\bigtriangleup \left| \nabla u\right| ^2
& = &\left| \nabla ^2u\right| ^2+\left\langle \nabla u,\nabla 
R\right\rangle +Ric\left( \nabla u,\nabla u\right) \nonumber \\
& \geq & \left| \nabla ^2u\right| ^2+\left\langle \nabla u,\nabla 
R\right\rangle \ .
\end{eqnarray}
For any $r>0$ and any $\overline{x}_0\in M$, by multiplying 
(\ref{6.18}) by the cutoff function in (\ref{2.7}) and integrating by
parts, we get
\begin{eqnarray}
\int_M\left| \nabla ^2u\right| ^2\varphi dx
& \leq & \frac 12\int_M\left| \nabla u\right| ^2\cdot \left| 
\bigtriangleup \varphi \right| dx+\int_M\left| \nabla ^2u\right| 
\cdot R\varphi dx \nonumber \\
&  & + \int_M\left| \nabla u\right| \cdot R\cdot \left| 
\nabla \varphi \right| dx \nonumber \\
& \leq & \frac{C_{16}^2}2\cdot \frac{C_3}{r^2}\int_M\varphi dx+
\frac 12\int_M\left| \nabla ^2u\right| ^2\varphi dx \nonumber \\
&  & +\frac 12\int_MR^2\varphi dx+C_{16}\cdot \left( 
\sup \limits_MR\right) \cdot \frac{C_3}r\int_M\varphi dx\ .\nonumber
\end{eqnarray}
Thus,
\begin{equation}
\label{6.19}
\int_M\left| \nabla ^2u\right| ^2\varphi dx\leq 
\left(C_3C_{16}^2\cdot \frac 1{r^2}+2C_{16}C_3\left( 
\sup \limits_MR\right) \frac
1r\right) \int_M\varphi dx+\int_MR^2\varphi dx\ . 
\end{equation}
By (\ref{6.1}), (\ref{2.7}) and the standard volume comparison, 
we have
\begin{eqnarray}
\label{6.20}
\int_MR^2\varphi dx&\leq&\left( \sup \limits_MR\right) 
\int_MR(x)e^{-\left( 1+\frac{d(x,\overline{x}_0)}r\right) }dx
\nonumber \\
& \leq & \left( \sup \limits_MR\right)\cdot \nonumber\\
&&\left( \int_{B(\overline{x}_0,r)}R(x)dx+\sum\limits_{k=0}^\infty 
e^{-2^{k-1}}\int_{B(\overline{x}_0,2^{k+1}r)\backslash B(\overline{x}_0,2^kr)}R(x)dx\right) \nonumber\\
&\leq&C_{17}r^2\log \left( 2+r\right)
\end{eqnarray}
and
\begin{eqnarray}
\label{6.21}
\int_M\varphi dx&\leq&\int_Me^{-\left( 1+\frac{d(x,\overline{x}_0)}r
\right) }dx \nonumber\\
& \leq & \int_{B(\overline{x}_0,r)}dx+\sum\limits_{k=0}^\infty 
e^{-2^{k-1}}\int_{B(\overline{x}_0,2^{k+1}r)\backslash 
B(\overline{x}_0,2^kr)}dx \nonumber\\
& \leq & C_{17}r^4
\end{eqnarray}
for some positive constant $C_{17}$ independent of $r$ 
and $\overline{x}_0.$

Substituting these two inequalities into (\ref{6.19}) we have 
\begin{equation}
\label{6.22}
\frac 1{r^4}\int_{B(\overline{x}_0,r)}\left| \nabla ^2u\right|
^2dx\leq C_{18}\left( \frac 1{r^2}+\frac 1r+\frac{\log (2+r)}{r^2}\right) 
\end{equation}
for some positive constant $C_{18}$ independent of $r$ and $\overline{x}_0.$

Since the holomorphic bisectional curvature of 
$g_{\alpha \overline{\beta }}$ is positive, it was shown 
in \cite{MSY} that the function 
$\left| \sqrt{-1}\partial \overline{\partial }u-\mbox{Ric}\right| ^2$ 
is subharmonic on $M$. Then by the mean value inequality and 
(\ref{6.20}), (\ref{6.22}), we have
\begin{eqnarray}
\left| \sqrt{-1}\partial \overline{\partial }u-Ric\right| ^2(\overline{x}_0)&\leq&\frac{C_{19}}{r^4}\int_{B(\overline{x}_0,r)}
\left| \sqrt{-1}\partial \overline{\partial }u-Ric\right| ^2(x)dx \nonumber\\
& \leq & \frac{2C_{19}}{r^4}\int_{B(\overline{x}_0,r)}\left( 
\left| \nabla ^2u\right| ^2+R^2\right) dx \nonumber\\
& \leq & C_{20}\left( \frac 1{r^2}+\frac 1r+\frac{\log (2+r)}{r^2}\right) \nonumber
\end{eqnarray}for some positive constants $C_{19},\ C_{20}$ independent of $r$
and $\overline{x}_0.$ Since $\overline{x}_0\in M$ and $r>0$ are arbitrary,
by letting $r\rightarrow +\infty $ we know that
$$
\sqrt{-1}\partial \overline{\partial }u= \mbox{Ric} \qquad 
\mbox{on} \quad M\ . 
$$
In summary, we have proved the following result.
\vskip 3mm{\bf \underline{Proposition 6.2}} \qquad 
Suppose $(M,g_{\alpha \overline{\beta }})$ is a complete noncompact
K\"ahler surface satisfying all the assumptions in the Main Theorem. 
Then there exists a strictly plurisubharmonic function $u(x)$ on $M$ 
satisfying the Poincar\'e--Lelong equation (\ref{6.11}) with the
estimate 
$$
|u(x)| \leq C\log\left(2+d(x,x_0)\right) \qquad \mbox{for all} \;\; x\in M
$$
for some positive constant C.

\section*{\S 7. Uniform estimates on multiplicity and the number 
of components of an ``algebraic'' divisor}

\setcounter{section}{7} \setcounter{equation}{0}

\qquad
Let $(M,g_{\alpha \overline{\beta }})$ be a complete noncompact
K\"ahler surface satisfying all the assumptions in the Main Theorem.
In this section, we will consider the algebra $P(M)$ of holomorphic 
functions of polynomial growth on $M$. We first construct 
$f_1,\ f_2$ in $P(M)$ which are algebraically independent over $\C.$

In the previous section, by solving the Poincar\'e--Lelong equation, 
we have obtained a strictly plurisubharmonic function $u$ 
on $M$ of logarithmic growth. As shown in \cite{Mo1}, the existence 
of nontrivial functions in the algebra $P(M)$ then follows readily 
from the $L^2$ estimates of the $\overline{\partial }$ operator on 
complete K\"ahler manifold of Andreotti--Vesentini \cite{AV}
and H\"ormander \cite{Ho}. For completeness, we give the proof as follows.

Let $x\in M$ and $\left\{ (z_1,z_2)\,|\ |z_1|^2+|z_2|^2<1\right\}$ 
be local holomorphic coordinates at $x$ with $z_1(x)=z_2(x)=0$. 
Let $\eta $ be a smooth cutoff function on $\C^2$ with 
$\mbox{Supp}\,\eta \subset \subset\left\{ |z_1|^2+|z_2|^2<1\right\}$ 
and $\eta \equiv 1$ on $\left\{|z_1|^2+|z_2|^2<\frac 14\right\} $. 
Then the function
$$
\eta \log |z|=\eta \left( z_1,z_2\right) \log 
\left(|z_1|^2+|z_2|^2\right) ^{\frac 12} 
$$
is globally defined on $M$ and is smooth except at $x$. Furthermore, 
the $(1,1)$ form $\partial \overline{\partial }\,(\eta \log |z|)$ is 
bounded from below. Since $u$ is strictly plurisubharmonic, we can 
choose a sufficiently large positive constant $C$ such that
$$
v=Cu+ 6 \eta \, \log \left|\, z\right|  
$$
is strictly plurisubharmonic on $M$. Then, for any nonzero tangent 
vector $\xi $ of type $(1,0)$ on $M$, we have
$$
\left\langle \sqrt{-1}\partial \overline{\partial }v+\mbox{Ric},\;
\xi \wedge \overline{\xi } \right\rangle >0\ . 
$$
Now $\overline{\partial }\,(\eta z_i)$, i=1,2 , is a $\overline{\partial }$
closed $(0,1)$ form on the complete K\"ahler manifold $M$. Using the
standard $L^2$--estimates of $\overline{\partial }$ operator (c.f. Theorem
2.1 in \cite{Mo1}), there exists a smooth function $u_i$ such that
$$
\overline{\partial }u_i=\overline{\partial }(\eta z_i)\ \qquad i=1,2 
$$
and
$$
\int_M\left| u_i\right| ^2e^{-v}dx\leq \frac 1c\int_M\left| \overline{%
\partial }(\eta z_i)\right| ^2e^{-v}dx 
$$
where $c$ is a positive constant satisfying 
$$
\left\langle \sqrt{-1} \partial \overline{\partial }v+ \mbox{Ric},\xi \wedge 
\overline{\xi }\right\rangle \geq c|\xi |^2
$$ 
whenever $\xi$ is a tangent vector on $\mbox{Supp}\, \eta $.
First of all, this estimate implies that $u_i$ is of polynomial growth
as the weight function $v$ is of logarithmic growth. Secondly, 
because of the singularity of $6\log \left| z\right| $ at $x$, it
forces the function $u_i$ and its first order
derivative to vanish at $x$. Therefore, the holomorphic functions 
$f_1=u_1-\eta z_1$ and $f_2=u_2-\eta z_2$ define a local biholomorphism 
at $x$. Clearly, they are algebraically independent over $\C$. 
This concludes our construction.

For later use, we also point out here that, as a consequence of the 
above argument, the algebra $P(M)$ separates points on $M$. In other
words, for any $x_1,x_2\in M$ with $x_1 \neq x_2$, 
there exists $f\in P(M)$ such that $f(x_1) \neq f(x_2).$

Before we can state our main result in this section, we need
the following definition.
For a holomorphic function $f\in P(M)$, we define the degree of $f$, 
$\deg (f)$, to be the infimum of all $q$ for which the following 
inequality holds
$$
\left| f(x)\right| \leq C(q)\left( 1+d^q(x,x_0)\right) \qquad 
\mbox{for all} \; x\in M, 
$$
where $x_0$ is some fixed point in $M$ and $C(q)$ is some positive 
constant depending on $q.$

Our main result in this section is the following uniform bound 
on the multiplicity of the zero divisor of a function $f\in P(M)$
by its degree.

\vskip 3mm{\bf \underline{Proposition 7.1}}
\qquad 
Let $(M,g_{\alpha \overline{\beta }})$ be a complete noncompact
K\"ahler surface as above. For $f\in P(M)$, let
$$
\left[ V\right] =\frac{\sqrt{-1}}{2\pi }\partial \overline{\partial }\log
\left| f\right| ^2 
$$
be the zero divisor, counting multiplicity, determined by $f$. Then, there
exists a positive constant $C$, independent of $f$, such that
$$
\mbox{mult} \left( \left[ V\right] ,x\right) \leq C\deg (f) 
$$
holds for all $x\in M.$

\vskip 3mm{\bf \underline{Proof.}} 
\qquad 
Recall that
the Ricci flow (\ref{2.1}) with $g_{\alpha \overline{\beta }}(x)$ as initial
metric has a solution $g_{\alpha \overline{\beta }}(x,t)$ for all times 
$t\in [0,+\infty )$ and satisfies the following estimates 
\begin{equation}
\label{7.1}
R(x,t) \leq \frac C{1+t}\  
\end{equation}
and 
\begin{equation}
\label{7.2}
\mbox{inj} \left( M,g_{\alpha \overline{\beta }}(\cdot ,t)\right) \geq
C_7(1+t)^{\frac 12}\  
\end{equation}
on $M\times [0,+\infty )$. 

Let $d_t$ be the distance function from an
arbitrary fixed point $\overline{x}_0\in M$ with respect to 
the metric $g_{\alpha \overline{\beta }}(\cdot ,t)$. 
By the standard Hessian comparison theorem (see \cite{ScY}), we have, 
for any unit real vector $v$ orthogonal
to the radial direction $\partial /\partial d_t$,
$$
\frac{\sqrt{\frac{\alpha _1}t}d_t}{\tan \left( \sqrt{\frac{\alpha _1}t}
d_t\right) }\leq \mbox{Hess}\left( d_t^2\right) (v,v)\leq \frac{\sqrt{\frac{\alpha
_1}t}d_t}{\tanh \left( \sqrt{\frac{\alpha _1}t}d_t\right) } \qquad
\mbox{when}\quad d_t\leq \frac \pi 4\sqrt{\frac t{\alpha _1}}\ . 
$$
Here, $\alpha _1$ is some positive constant depending only on the 
constants $C $ and $C_7$ in (\ref{7.1}) and (\ref{7.2}). 
Hence, for any unit vector $\widetilde{v}$, we have
$$
\frac{\sqrt{\frac{\alpha _1}t}d_t}{\tan \left( \sqrt{\frac{\alpha _1}t}
d_t\right) }\leq \mbox{Hess}\left( d_t^2\right) (\widetilde{v},\widetilde{v}
) + \mbox{Hess} \left( d_t^2\right) (J\widetilde{v},J\widetilde{v})
\leq \frac{2\sqrt{\frac{\alpha _1}t}d_t}{\tanh \left( 
\sqrt{\frac{\alpha _1}t}d_t\right) } 
$$
whenever $d_t\leq \frac \pi 4\sqrt{\frac t{\alpha _1}}$. Since $M$ is
K\"ahler, the above expression is equivalent to
$$
\frac{\sqrt{\frac{\alpha _1}t}d_t}{\tan \left( \sqrt{\frac{\alpha _1}t}%
d_t\right) }\omega _t\leq \sqrt{-1}\partial \overline{\partial }d_t^2\leq 
\frac{2\sqrt{\frac{\alpha _1}t}d_t}{\tanh \left( \sqrt{\frac{\alpha _1}t}%
d_t\right) }\omega _t\ . 
$$
In particular, we have 
\begin{equation}
\label{7.3}
\frac 12\omega _t\leq \sqrt{-1}\partial \overline{\partial }
d_t^2\leq 4\omega _t \qquad \mbox{whenever} \quad d_t \leq \frac {\pi} {4}
\sqrt{\frac t{\alpha _1}}\ . 
\end{equation}
Here, $\omega _t$ is the K\"ahler form of the metric 
$g_{\alpha \overline{\beta }}(\cdot ,t)$.

We next claim that 
\footnote{we are grateful to Professor L.F.Tam for this suggestion} 
\begin{equation}
\label{7.4}\sqrt{-1}\partial \overline{\partial }\log \tan \left( \sqrt{
\frac{\alpha _1}t}\frac{d_t}2\right) \geq 0\ ,\qquad 
\mbox{whenever} \quad d_t\leq
\frac \pi 4\sqrt{\frac t{\alpha _1}}\ . 
\end{equation}

In fact, after recaling, we may assume that the sectional curvature of 
$g_{\alpha \overline{\beta }}(\cdot ,t)$ is less than $1$ and 
$\sqrt{\frac{\alpha _1}t}=1$. Then by the standard Hessian comparison, 
we have
$$
\mbox{Hess} \left( d_t\right) (v,v)\geq \frac 1{\tan d_t}\left( 
\left| v\right|_t^2-\left\langle v,\frac \partial 
{\partial d_t}\right\rangle _t^2\right) 
$$
for any vector $v$ and $d_t\leq \frac \pi 4$. Thus, by a direct computation,
\begin{eqnarray}
&  & \mbox{Hess}\left( \log \tan \left( \frac{d_t}2\right) \right) (v,v)+
\mbox{Hess} \left( \log \tan \left( \frac{d_t}2\right) \right) (Jv,Jv) \nonumber\\
& \geq &\frac 1{\left( \tan d_t\right) \tan \left( \frac{d_t}2\right) }
\left( 1+\tan d_t\right) \left| v\right| _t^2 \nonumber\\
& \geq & 0\ ,\nonumber
\end{eqnarray}
which is our claim (\ref{7.4}).

Now for any $0<b<a<\frac \pi 8\sqrt{\frac t{\alpha _1}}$, it follows from
Stoke's theorem that
\begin{eqnarray}
0&\leq&\sqrt{-1}\int\limits_{\left\{ b\leq d_t\leq a\right\} }\left[ V\right] \wedge\partial\overline\partial \log \tan \left( \sqrt{\frac{\alpha _1}t}\frac{d_t}2\right) \nonumber\\
& = &\sqrt{-1}\int\limits_{\left\{ d_t=a\right\} \ }\left[ V\right] 
\wedge \overline{\partial }\log \tan \left( \sqrt{\frac{\alpha _1}t}\frac{d_t}2\right) \nonumber\\
&&-\sqrt{-1}\int\limits_{\left\{ d_t=b\right\} \ }\left[ V\right] 
\wedge \overline{\partial }\log \tan \left( \sqrt{\frac{\alpha _1}t}
\frac{d_t}2\right) \ .\nonumber
\end{eqnarray}
Then, it is no hard to see that for $0<b<a<\frac \pi 8\sqrt{\frac t{\alpha _1}}$,  
\begin{equation}
\label{7.5}
\frac{\sqrt{-1}}{a^2}\int\limits_{\left\{ d_t=a\right\} \ \
}\left[ V\right] \wedge \overline{\partial }\left( d_t^2\right) \geq \frac
12\cdot \frac{\sqrt{-1}}{b^2}\int\limits_{\left\{ d_t=b\right\} \ \ }\left[
V\right] \wedge \overline{\partial }\left( d_t^2\right) \ . 
\end{equation}

Using Stoke's theorem on the right hand side of (\ref{7.5}) and letting 
$b\rightarrow 0$ it follows from the inequality of Bishop--Lelong that 
\begin{equation}
\label{7.6}
\frac{\sqrt{-1}}{a^2}\int\limits_{\left\{ d_t=a\right\} \ \
}\left[ V\right] \wedge \overline{\partial }\left( d_t^2\right) \geq \alpha
_2\mbox{mult}\,\left( \left[ V\right] ,\overline{x}_0\right) \ , 
\end{equation}
for some positive absolute constant $\alpha _2$.

Then by (\ref{7.3}), (\ref{7.5}), (\ref{7.6}) and Stoke's theorem, we have%
\begin{eqnarray}
\label{7.7}
&   &\frac 1{a^2}\int\nolimits_{B_t(\overline{x}_0,a)\backslash 
B_t(\overline{x}_0,\frac a2)}\left[ V\right] \wedge \omega _t\nonumber\\
& \geq &\frac 1{4a^2}\int\nolimits_{B_t(\overline{x}_0,a)\backslash 
B_t(\overline{x}_0,\frac a2)}\left[ V\right] \wedge \sqrt{-1}\partial 
\overline{\partial }\left( d_t^2\right) \nonumber\\
& = &\frac{\sqrt{-1}}4\left( \frac 1{a^2}\int\nolimits_{\{d_t=a\}}
\left[ V\right] \wedge \overline{\partial }\left( d_t^2\right) -
\frac 1{4\cdot \left( \frac a2\right) ^2}\int\nolimits_{\{d_t=\frac a2\}}
\left[ V\right] \wedge \overline{\partial }
\left( d_t^2\right) \right)\nonumber\\
& = &\frac{\sqrt{-1}}8\left( \frac 1{a^2}\int\nolimits_{\{d_t=a\}}
\left[ V\right] \wedge \overline{\partial }\left( d_t^2\right) -
\frac 1{2\cdot \left( \frac a2\right) ^2}\int\nolimits_{\{d_t=\frac a2\}}
\left[ V\right] \wedge \overline{\partial }\left( d_t^2\right) \right)\nonumber\\
&    &+\frac{\sqrt{-1}}{8a^2}\int\nolimits_{\{d_t=a\}}\left[ V\right] 
\wedge \overline{\partial }\left( d_t^2\right)\nonumber\\ 
& \geq &\frac{\sqrt{-1}}{8a^2}\int\nolimits_{\{d_t=a\}}\left[ V\right] 
\wedge \overline{\partial }\left( d_t^2\right)\nonumber\\
& \geq &\frac{\alpha _2}8\mbox{mult}\,\left( \left[ V\right] ,\overline{x}_0\right)
\end{eqnarray}
for $0<a<\frac \pi 8\sqrt{\frac t{\alpha _1}}.$

For the function $f\in P(M)$, let $\widetilde{x}_0$ be a point close to 
$\overline{x}_0$ such that $f(\widetilde{x}_0)\neq 0$. By definition, for any 
$\delta >0$, there exists a constant $C(\delta )>0$ such that 
\begin{equation}
\label{7.8}\left| f(x)\right| \leq C(\delta )\left( 1+d_0^{\deg (f)+\delta
}(x,\widetilde{x}_0)\right) \qquad \mbox{on} \quad M\ . 
\end{equation}
By equation (\ref{2.1}) and estimate (\ref{7.1}), we have
\begin{eqnarray}
\frac{\partial g_{\alpha \overline{\beta }}(\cdot ,t)}
{\partial t}&\geq&-R(\cdot ,t)g_{\alpha \overline{\beta }}(\cdot ,t)
\nonumber \\
& \geq &-\frac C{1+t}g_{\alpha \overline{\beta }}(\cdot ,t)\ ,\nonumber
\end{eqnarray}which implies that
$$
g_{\alpha \overline{\beta }}(\cdot ,0)\leq (1+t)^Cg_{\alpha \overline{\beta }%
}(\cdot ,t)\qquad \mbox{for any} \quad t>0\ . 
$$
Hence, (\ref{7.8}) becomes 
\begin{equation}
\label{7.9}
\left| f(x)\right| \leq C(\delta )\left\{ 1+\left[ (1+t)^{\frac
C2}d_t(x,\widetilde{x}_0)\right] ^{\deg (f)+\delta }\right\} \qquad
\mbox{on} \quad M\ . 
\end{equation}

We now fix $t=\frac{\alpha _1}{\pi ^2}4^{K+8}$ for each positive 
interger $K$. Set
$$
v_K(x)=\int\nolimits_{B_t(\widetilde{x}_0,2^K)}-G_t^{(K)}(x,y)\bigtriangleup
_t\log \left| f(y)\right| ^2\cdot \omega _t^2(y)\ , 
$$
where $G_t^{(K)}$ is the positive Green function with value zero on the
boundary $\partial B_t(\widetilde{x}_0,2^K)$ with respect to the metric 
$g_{\alpha \overline{\beta }}(\cdot ,t)$. The function 
$\log \left| f\right|^2-v_K$ is then harmonic on $B_t(\widetilde{x}_0,2^K)$. 
From the maximum principle and (\ref{7.9}), we have
\begin{eqnarray}
\label{7.10}
\log \left( \left| f(\widetilde{x}_0)\right| ^2\right) -
v_K(\widetilde{x}_0)&\leq&\sup \limits_{x\in \partial B_t
(\widetilde{x}_0,2^K)}\log \left| f(x)\right| ^2    \nonumber\\
&\leq&C_{19}K\left( \deg (f)+\delta \right) +C^{\prime }(\delta )
\end{eqnarray}
for some positive constants $C_{19}$, $C^{\prime }(\delta )$
independent of $K,\ f,$ and $\overline{x}_0.$

On the other hand, since the volume growth condition (i) is preserved for
all times, by virtue of (\ref{6.4}) 
(c.f. Proposition 1.1 in \cite{Mo1}), we have
\begin{eqnarray}
-v_K(\widetilde{x}_0)
& \geq & \frac 1{C_9}\int\nolimits_{B_t(\widetilde{x}_0,
          2^K)}\frac 1{d_t^2(x,\widetilde{x}_0)}
          \bigtriangleup _t\log \left| f(x)\right| ^2\cdot \omega _t^2(x)
          \nonumber\\
& \geq &  \frac 1{C_9}\sum\limits_{j=1}^K\left( \frac 1{2^j}
      \right)^2\int\nolimits_{B_t(\widetilde{x}_0,2^j)\backslash 
          B_t(\widetilde{x}_0,2^{j-1})}\bigtriangleup _t
           \log \left| f(x)\right| ^2\cdot \omega _t^2(x)\ . \nonumber
\end{eqnarray}
Then, by (\ref{7.7}) and the fact that $\widetilde{x}_0$ is
arbitrarily close to $\overline{x}_0,$%
\begin{equation}
\label{7.11} 
-v_K(\widetilde{x}_0)\geq C_{20}\,K\,\mbox{mult}\,\left( 
\left[ V\right] , \overline{x}_0\right) 
\end{equation}
for some positive constant $C_{20}$ independent of $K,\ f$ and 
$\overline{x}_0.$

Therefore, by combining (\ref{7.10}) and (\ref{7.11}) and letting 
$K\rightarrow +\infty $ and then $\delta \rightarrow 0$, we obtain
$$
\mbox{mult}\,\left( \left[ V\right] ,\overline{x}_0\right) 
\leq C_{21}\deg (f) 
$$
where $C_{21}$ is some positive constant independent of $f$ 
and $\overline{x}_0.$
\hfill $\Box $

\vskip 3mm 
A modified version of the proof of Proposition 7.1 gives the 
uniform bound on the number of irreducible components of $[V].$

\vskip 3mm {\bf \underline{Proposition 7.2}} 
\qquad 
Suppose $(M,g_{\alpha \overline{\beta }})$ is a complete noncompact 
K\"ahler surface as assumed in Proposition 7.1. Let $f$ be a holomorphic 
function of polynomial growth, 
$$\left[ V\right] =\frac{\sqrt{-1}}{2\pi }
\partial \overline{\partial }\log\left| f\right| ^2
$$ 
be the corresponding zero divisor determined by $f$.
Then the number of irreducible components of $[V]$ is not bigger 
than $C\deg(f)$ for the same positive constant $C$ as in Proposition 7.1.

\vskip 3mm{\bf \underline{Proof.}} 
\qquad 
Let $g_{\alpha \overline{\beta }}(\cdot ,t)$ be
the evolving metric to the Ricci flow with 
$g_{\alpha \overline{\beta }}(\cdot)$ 
as the initial metric and $[V_1],\ [V_2],\ \cdots ,\ [V_l]$ be any $l$
distinct irreducible components of $\left[ V\right] $. Fix a
constant $a>0$ such that the intersection of the smooth points of 
$\left[ V_i \right]$ with $B_0(\overline{x}_0,a)$ is nonempty for each 
$0 \leq i \leq l.$

Choose $t=\frac{\alpha _1}{\pi ^2}4^{K+8}a^2$ for each positive 
integer $K$. As the manifold $M$ is Stein by Theorem 5.1, each 
$\left[ V_i\right] $ must be noncompact. Hence, for $j=1,2,\cdots ,K$, 
we have
$$
\left[ V_i\right] \cap \left( B_t(\overline{x}_0,2^ja)\left\backslash 
B_t( \overline{x}_0,2^{j-1}a)\right. \right) \neq \emptyset 
$$
and there exists a point $x_j\in \left[ V_i\right] $ with 
$d_t(x_j,\overline{x}_0)=\frac 322^{j-1}a$ in the middle of $B_t(\overline{x}_0,2^ja)\left\backslash B_t( 
\overline{x}_0,2^{j-1}a)\right.$.
%By perturbing $x_j$ a little bit if necessary, we may assume
%that $x_j$ is a smooth point on $\left[ V_j\right]$.
The triangle inequality says 
$$
B_t(x_j, 2^{j-2}a)\subset \left( B_t(\overline{x}_0,2^ja)\left\backslash
B_t(\overline{x}_0,2^{j-1}a)\right. \right) \ . 
$$
Applying a slight variant of (\ref{7.7}) to $\left[ V_i\right] $, we have
\begin{eqnarray}
\frac 1{\left( 2^{j-2}a\right) ^2}\int\nolimits_
{B_t(x_j, 2^{j-2}a)}\left[ V_i\right] \wedge \omega _t
& \geq & \frac{\alpha _2}8 \mbox{mult}\,\left( 
\left[ V_i\right] ,x_j\right)  \nonumber\\
& \geq & \frac{\alpha _2}8\ . \nonumber
\end{eqnarray}
%as $x_j$ is a smooth point. 
Since $\sum\limits_{i=1}^l\left[ V_i\right] $ 
is only a part of the divisor $\left[ V\right] $, we get
$$
\frac 1{\left(2^{j-2}a\right) ^2}\int\nolimits_{B_t(\overline{x}
_0,2^ja)\left\backslash B_t(\overline{x}_0,2^{j-1}a)\right. }
\bigtriangleup_t\log \left| f(x)\right| ^2\cdot \omega _t^2(x)
\geq \frac{\alpha _2}8\, l\ . 
$$

The subsequent argument is then exactly as in the proof of 
Proposition 7.1. In the end, we have
$$
C_{20}K\cdot l\leq -\log \left( \left| f\left( \widetilde{x}_0\right)
\right| ^2\right) +C_{19}K\left( \deg (f)+\delta \right) +
C^{\prime }(\delta)\ . 
$$
Letting $K\rightarrow +\infty $ and then $\delta \rightarrow 0$, we get
the desired estimate.
\hfill$\Box $

\section*{\S8. Proof of the main theorem}

\setcounter{section}{8} \setcounter{equation}{0}

\qquad 
In this section, we will basically follow the approach of 
Mok in \cite{Mo1}, \cite{Mo3} to accomplish the proof of the 
main theorem. Let $M$ be a K\"ahler surface as assumed in the 
Main Theorem. Recall that $P(M)$ stands for the algebra of 
holomorphic functions of polynomial growth on $M$. Let $R(M)$ be 
the quotient field of $P(M)$. By an abuse of terminology, we will
call it the field of rational functions on $M$.

In the previous section, we showed that there exist two functions 
$f_1,\ f_2\in P(M)$ giving local holomorphic coordinates at any 
given point $x \in M$, and that the algebra $P(M)$ separates points
on $M$. Moreover, we obtained the following basic multiplicity estimate 
\begin{equation}
\label{8.1}
\mbox{mult}\,\left( \left[ V\right] ,x\right) \leq C\deg (f)
\end{equation}
for all $x\in M$ and $f\in P(M)$, where 
$$
\left[ V\right] =\frac{\sqrt{-1}}{2\pi }\partial 
\overline{\partial }\log \left| f\right| ^2
$$ 
is the zero divisor of $f$ and $C$ is a constant independent of
$f$ and $x$. Thus, by combining these facts with 
the classical arguments of Poincar\'e and Siegel, we have (c.f. 
the proof of Proposition 5.1 in \cite{Mo1}) 
\begin{equation}
\label{8.2}
\dim {}_{{\C}}H_p\leq 10^3Cp^2\ ,
\end{equation}
where $H_p$ denotes the vector space of all holomorphic functions 
with degree $\leq p$, and the field of rational functions $R(M)$ 
is a finite extension field over ${\C}(f_1,f_2)$ for some algebraically 
independent holomorphic functions $f_1,\ f_2\in P(M)$ over $\C$.
By the primitive element theorem, we can then write
$$
R(M)={\C}\left( f_1,f_2,\frac{f_3}{f_4}\right) 
$$
for some $f_3,\ f_4\in P(M)$. 

Now, consider the mapping $F:M\rightarrow {\C}^4 $ defined by
$$
F=\left( f_1,f_2,f_3,f_4\right) \ . 
$$
Since $R(M)$ is a finite extension field of ${\C}(f_1\ f_2)$, 
$f_3$ and $f_4$ satisfy equations of the form
$$
f_3^p+\sum\limits_{j=0}^{p-1}P_j(f_1,f_2)f_3^j=0\ , 
$$
$$
f_4^q+\sum\limits_{j=0}^{q-1}Q_j(f_1,f_2)f_4^j=0\ , 
$$
where $P_j(w_1,w_2)$, $Q_j(w_1,w_2)$ are rational functions of 
$w_1,\ w_2$. After clearing denominators, we see that
$f_1, f_2, f_3, f_4$ satisfy polynomial equations
$$
P(f_1,f_2,f_3,f_4) = 0 \qquad \mbox{and} \qquad 
Q(f_1,f_2,f_3,f_4) = 0.
$$ 
Let $Z_0$ be the subvariety of $\C^4$ defined by
$$
Z_0=\left\{ \left( w_1,w_2,w_3,w_4\right)\in \C^4 \left| \ 
\begin{array}{l}
\bigbreak P(w_1,w_2,w_3,w_4) = 0\  \\ 
Q(w_1,w_2,w_3,w_4) = 0 
\end{array}
\right. \right\} \ , 
$$
and let $Z$ be the connected component of $Z_0$ containing $F(M)$. 
It is clear that $\dim {}_{{\C}}Z=2.$

In the following we will show that $F$ is an "almost injective"
and "almost surjective" map to $Z$ and we can desingularize $F$ to 
obtain a biholomorphic map from $M$ onto a quasi--affine algebraic 
variety by adjoining a finite number of holomorphic functions of 
polynomial growth.

First of all, we claim that $Z$ is irreducible and $F$ is 
"almost injective", i.e., there exists a subvariety $V$ of $M$ such that 
$F|_{M\backslash V}:M\backslash V\rightarrow Z$ is an
injective locally biholomorphic mapping.
Indeed, take $V$ to be the union of $F^{-1}(\mbox{Sing}(Z))$ and the 
branching locus of $F$, here Sing(Z) denotes the singular set
of Z. It is clear that $F$ is locally biholomorphic on 
$M\backslash V$. That $F$ is also injective there follows from 
the fact that $P(M)$ separates points and $f_1,...,f_4$ generate
$P(M)$. To see the irreducibility of $Z$, note that 
$M \backslash F^{-1}(\mbox{Sing}(Z))$
is connected and hence $\overline{F(M \backslash F^{-1}(\mbox{Sing}(Z))}$ 
is irreducible (as its set of smooth points is connected). Since 
$F(M) \subset \overline{F(M \backslash F^{-1}(\mbox{Sing}(Z))}$, by
the definition of $Z$, it must be irreducible.
%$F(x_1)\neq F(x_2)$ for any different two points $x_{1,}\ x_2$. Let $V$ be
%the union of the branching locus of $F$ and $F^{-1}(sing(Z))$. Then $F$ is
%injective and locally biholomorphic on $M\backslash V$. By codimension
%reason, $\overline{F(M\backslash V)}$ is irreducible ( i.e., its smooth
%points are connected ). This says, $F(M)$ lies in an irreducible component,
%so $Z$ is irreducible by definition.

Next, we come to the "almost surjectivity" of $F$, i.e., there 
exists an algebraic subvariety $T$ of $Z$ such that $F(M)$ contains 
$Z\backslash T$. The method of Mok\cite{Mo1} in proving the
almost surjectivity of $F$ is to solve an ideal 
problem for each $x \in Z \backslash T_0$ missed by $F$, where 
$T_0$ is some fixed algebraic subvariety of $Z$ containing the singular set 
of $Z$. The solution of the ideal problem gives a holomorphic function
$f_{x}\in P(M)$ with degree bounded independent of $x$ which will 
correspond to a rational function on $\C^4$ with pole set passes through 
$x$. Then, the almost surjectivity of $F$ follows.  
Otherwise, one could select an infinite number of linearly independent 
$f_x$'s contradicting the finite dimensionality of the space of
holomorphic functions with polynomial growth of some fixed 
degree, c.f. (\ref{8.2}).

In \cite{Mo1}, Mok used the solution $u$ of the Poincar\'e--Lelong 
equation as the weight function in the Skoda's estimates for solving 
the ideal problem. In his case, because of his curvature quadratic
decay condition, the growth of $u$ is bounded both from above
and below by the logarithm of the distance function on $M$. This
does not work in our case because we do not have the luxury of the
lower bound of $u$. However, thanks to the Steinness
of $M$ by Theorem 5.1, we can adapt the argument of Mok in \cite{Mo3} 
to choose another weight function by resorting to Oka's theory of pseudoconvex 
Riemann domains.

Before we carry out the above procedures in proving the almost
surjectivity of $F$. We first need
%Before we start to show that $F$ ``almost surjective'', we would like 
to construct a nontrivial holomorphic $(2,0)$ vector field of polynomial
growth on $M$. 

Consider the anticanonical line bundle, ${\bf K} ^{-1}$, on 
$M$ equipped with the induced Hermitian metric, its curvature 
form $\Omega ({\bf K} ^{-1})$ is then simply the Ricci form of $M$. 
Let $u$ be the strictly plurisubharmonic function of logarithmic growth 
obtained in Proposition 6.2. For any given point $\overline{x}_0\in M$, 
let $\left\{ z_1,z_2\right\} $ be local holomorphic coordinates at 
$\overline{x}_0$. Choose a smooth cutoff function $\eta $ supporting
in this local holomorphic coordinate chart with value 
$1$ in a neighborhood of $\overline{x}_0$.
We study the following $\overline{\partial }$ 
equation for the sections of ${\bf K} ^{-1}$ on $M$,
\begin{equation}
\label{8.3}
\overline{\partial }S=\overline{\partial }\left( \eta \frac
\partial {\partial z_1}\wedge \frac \partial {\partial z_2}\right).
\end{equation}
Clearly, we can choose $k>0$ large enough such that
$$
k\sqrt{-1}\partial \overline{\partial }u+\Omega({\bf K} ^{-1})+
3\sqrt{-1}\partial \overline{\partial }\left( \eta \log 
\left( |z_1|^2+|z_2|^2\right) \right)>0\ . 
$$
Then by the standard $L^2$ estimate of $\overline{\partial }$ 
operator on Hermitian holomorphic line bundles (c.f. Theorem 1.2 
in \cite{Mo3}), equation (\ref{8.3}) has a smooth solution $S(x)$ 
satisfying the estimate
\begin{eqnarray}
\label{8.4}
& &\int_M\left| S\right| ^2e^{-ku-3\eta \log \left( 
|z_1|^2+|z_2|^2\right) }\omega ^2    \nonumber\\
& \leq & C\int_M\left| \overline{\partial }\left( 
\eta \frac \partial {\partial z_1}\wedge \frac \partial 
{\partial z_2}\right) \right| ^2e^{-ku-3\eta \log \left( 
|z_1|^2+|z_2|^2\right) }\omega ^2\nonumber\\
& < & + \infty
\end{eqnarray}
for some positive constant $C$. Recall the Poincar\'e--Lelong equation 
for the section $S(x)$ of the anticanonical line bundle $\bf{K}^{-1}$,
$$
\frac{\sqrt{-1}}{2\pi }\partial \overline{\partial }\log \left| S\right|
^2=[V]-\frac 1{2\pi }\mbox{Ric} \qquad \mbox{on} \quad M\ , 
$$
where $\left[ V\right] $ is the zero divisor of $S(x)$ (c.f. \cite{Mo3}).
Thus, $\log \left| S\right| ^2+u$ is subharmonic and so is 
$\left| S\right|^2e^u=\exp (\log \left| S\right| ^2+u).$ 
Since $M$ has positive Ricci
curvature and maximal volume growth, we can apply the mean value
inequality of subharmonic functions, (\ref{8.4}) and the fact that $u$ has
logarithmic growth to show that $S(x)$ is of polynomial growth. Set 
$$
v=\eta \left( \frac \partial {\partial z_1}\wedge \frac \partial
{\partial z_2}\right) -S\ .
$$
Then, $v$ is a nontrivial holomorphic $(2,0)$ vector field over $M$ with
polynomial growth we desired.

%According to Mok \cite{Mo1} to show that the mapping $F$ is ``almost
%surjective'', i.e., there exists an algebraic subvariety $T$ of $Z$ such
%that $F(M)$ contains $Z\backslash T$, one only need to solve an ideal
%problem for each $x\in Z$ missed $F$, except for a certain algebraic
%subvariety $T$, of $Z$, containing the singularities of $Z$. More precisely,
%one only need to prove that for each $x\notin F(M)\cup T$, there exists some 
%$f_x\in P(M)$ of degree bounded independent of $x$, which is the pull--back
%under $F$ of some rational function whose pole set pass through $x$. In fact
%if this was done, then the mapping $F$ must be ``almost surjective''
%otherwise one could select an infinite family of linear indendent of $f_x$%
%'s, giving a contradiction with (\ref{8.2}).
%In \cite{Mo1}, Mok used the solution $u$ of the Poincar\'e--Lelong equation
%as the weight function in the Skoda's estimates for solving the ideal
%problem, where the facts that the function $u$ is bounded above by logarithm
%and is also bounded below positively by logarithm are necessary for the
%argument. But in our situation we don't know whether the function $u$ has a
%positive lower bound of logarithm growth. To resolve the difficulty we now
%adapt the argument of Mok in \cite{Mo3} to choose another weight function by
%resorting to Oka's theory of pseudoconvex Riemann domains.

Now, for any $f_i,f_j\in \{f_1,f_2,f_3,f_4\}$ with $df_i\wedge
df_j\not \equiv 0$, we can choose the point $\overline{x}_0$ in the
above construction of $v$ so that the holomorphic function $f_{ij}$ defined
by
\begin{equation}
\label{8.5}
f_{ij}=\left\langle v,df_i\wedge df_j\right\rangle  
\end{equation}
is a nontrivial holomorphic function of polynomial growth. Here, 
we have used the fact that $\left\| df_i\wedge df_j\right\| $ 
grows at most polynomially by the gradient estimate of harmonic functions 
of Yau \cite{Y1}. It is obvious that the zero divisor of 
$df_i\wedge df_j$ is contained in the zero divisor of $f_{ij}$, for
which we denote by $V_0$. Since $M$ is Stein, the same is also
true for $M \backslash V_0$.

%Suppose $dw_i\wedge dw_j\not \equiv 0$ on $Z\,(\subset {\C^4})$, for some $%
%1\leq i,j\leq 4$. Recall that $M$ is a Stein manifold. Let $V_0$ be the zero
%divisor of the above nontrivial holomorphic functions $f_{ij}=\left\langle
%v,df_i\wedge df_j\right\rangle $ of polynomial growth. Thus $M\backslash V_0$
%is also a Stein manifold. 

Denote by $\pi_{ij}:Z\rightarrow \C^2$ the projection map given by 
$(w_1,w_2,w_3,w_4) \mapsto (w_i,w_j)$. Then, the map
$$
\rho =\pi _{ij}\circ F:M\backslash V_0\rightarrow \C^2 
$$
realises the Stein manifold $M\backslash V_0$ as a Riemann domain of
holomorphy over $\C^2.$
%On the Riemann domain of holomorphy $\rho :M\backslash V_0\rightarrow \C^2$,
Let $\delta (x)$ be the Euclidean distance to the boundary as in Oka 
\cite{O}. Then, $-\log \delta $ is a plurisubharmonic function on 
$M\backslash V_0$ by a theorem of Oka \cite{O}. $\delta (x)$ will be
used in the weight function of the Skoda's estimate mentioned 
above. It is essential to estimate it from below in terms of the
intrinsic distance $d(x,x_0)$ on $M$.

%that $\delta (x)$ can be estimated in terms of the intrinsic distance 
%$d(x,x_0)$ on $M$ and the holomorphic function $f_{ij}$. 

\vskip 3mm{\bf \underline{Lemma 8.1}} 
\qquad 
There exist positive constants $p$ and $C$ such that
$$
\delta (x)\geq C\left| f_{ij}(x)\right| ^2\left( d(x,x_0)+1\right) ^{-p}\ . 
$$

\vskip 3mm{\bf \underline{Proof.}} 
\qquad 
Let $v_i,\ v_j$ be two holomorphic vector fields on $M \backslash V_0$ defined 
by
$$
\left\langle v_k,df_l\right\rangle =\delta _{kl}\ ,\qquad k,l=i,j\ . 
$$
By the Cramer's rule, we have
$$
\left| v_k\right| \leq \frac{\left| df_i\right| +\left| 
df_j\right| }{\left|df_i\wedge df_j\right| }\leq 
\frac{\left| v\right| \left( \left| df_i\right|+\left| 
df_j\right| \right) }{\left| f_{ij}\right| }\leq C_{22}
\frac{\left(d(x,x_0)+1\right) ^{k_1}}{\left| f_{ij}(x)\right| }
\quad \mbox{on}\quad M\backslash V_0, 
$$
for $k=i,j$ and some positive constants $C_{22}$, $k_1$.

Since $f_{ij}$ is of polynomial growth, $|\nabla f_{ij}|$ is also of
polynomial growth by the gradient estimate of Yau, i.e.,%
$$
\max \left\{ f_{ij}(x),\left| \nabla f_{ij}(x)\right| \right\} \leq
C_{23}\left( d(x,x_0)+1\right) ^{k_2} \qquad \mbox{on} \quad M\ , 
$$
for some positive constants $C_{23}$ and $k_2$. 
Take $x \in M\backslash V_0$, then for any 
$y\in B\left( x,\left| f_{ij}(x)\right| \left/ 3C_{23}\left(
d(x,x_0)+1\right) ^{k_2}\right. \right) $, we have
\begin{eqnarray}
\label{8.6}
\left| f_{ij}(y)\right|&\geq&\left| f_{ij}(x)\right| -C_{23}
\left( d(x,x_0)+2\right) ^{k_2}\cdot \frac{\left| f_{ij}(x)
\right| }{3C_{23}\left( d(x,x_0)+1\right) ^{k_2}}
\nonumber\\
& \geq & \frac 12\left| f_{ij}(x)\right| \ .
\end{eqnarray}
This implies
$$
B\left( x,\left| f_{ij}(x)\right| \left/ 3C_{23}\left(
d(x,x_0)+1\right) ^{k_2}\right. \right) \subset M\backslash V_0
$$
and
\begin{equation}
\label{8.7}
\left| v_k(y)\right| \leq 2C_{22}\frac{\left( d(x,x_0)+1\right)
^{k_1}}{\left| f_{ij}(x)\right| }\ , 
\end{equation}
for all 
$y\in B\left( x,\left| f_{ij}(x)\right| \left/ 3C_{23}\left(
d(x,x_0)+1\right) ^{k_2}\right. \right) $, $k=i,j$.

By the definition of $\delta (x)$, it suffices to prove
\begin{eqnarray}
\label{8.8}
& &\rho \left( B\left( x,\left| f_{ij}(x)\right| \left/ 6C_{23}
\left( d(x,x_0)+1\right) ^{k_2}\right. \right) \right) \nonumber\\
& \supset & B_{\C^2}\left( \rho(x),C_{24}\left| f_{ij}(x)\right|^2 
\left/ \left( d(x,x_0)+1\right) ^{k_1+k_2}\right. \right),
\end{eqnarray}
for some positive constant $C_{24}$. Here, $B_{{\C^2}}(a,r)$
denotes the Euclidean ball in $\C^2$ with center $a$ and radius $r$.

Consider the real vector field
\begin{eqnarray}
\xi & = & \alpha _i\left( 2\mbox{Re}\,\left( v_i\right) \right) +
\alpha _j\left( 2\mbox{Re}\,\left( v_j\right) \right) +
\beta _i\left( 2\mbox{Im}\,\left( v_i\right) \right) +
\beta _j\left( 2\mbox{Im}\,\left( v_j\right) \right)  \nonumber\\
& = & \left( \alpha _i-\sqrt{-1}\beta _i\right) v_i+
\left( \alpha _j-\sqrt{-1}\beta _j\right) v_j+
\left( \alpha _i+\sqrt{-1}\beta _i\right) \overline{v}_i \nonumber\\
& & +\left( \alpha _j+\sqrt{-1}\beta _j\right) \overline{v}_j \nonumber
\end{eqnarray}
with $\left| \alpha _i\right| ^2+\left| \alpha _j\right| ^2+
\left| \beta_i\right| ^2+\left| \beta _j\right| ^2=1$. 
Clearly $\xi $ also satisfies (\ref{8.7}). Let 
$\gamma _\xi (\tau )$ be the integral curve in $M$ defined by 
$\xi$ and passes through $x$, i.e.,
\begin{equation}
\label{8.9}
\left\{ 
\begin{array}{l}
\bigbreak \displaystyle \frac{d\gamma _\xi (\tau )}{d\tau }=\xi \\ \gamma
_\xi (0)=x\ . 
\end{array}
\right. 
\end{equation}
We have
$$
\frac{d\left( f_i\circ \gamma _\xi (\tau )\right) }
{d\tau }=\left\langle \xi,df_i\right\rangle =
\alpha _i-\sqrt{-1}\beta _i\ , 
$$
$$
\frac{d\left( f_j\circ \gamma _\xi (\tau )\right) }
{d\tau }=\left\langle \xi,df_j\right\rangle =
\alpha _j-\sqrt{-1}\beta _j\ , 
$$
and  
\begin{equation}
\label{8.10}
\left| f_i\circ \gamma _\xi (\tau )-f_i(x)\right| ^2+\left|
f_j\circ \gamma _\xi (\tau )-f_j(x)\right| ^2=\tau ^2\ . 
\end{equation}
Note that (\ref{8.10})\ implies that $\gamma _\xi (\tau )$ cannot 
always stay in 
$$
B\left(x,\left|f_{ij}(x)\right|\left/6C_{23}\left(d(x,x_0)+1
\right) ^{k_2}\right.\right),
$$ 
otherwise $F=(f_1,f_2,f_3,f_4)$ would become unbounded in this ball. 
Denote by $\tau_0$ the first time when $\gamma _\xi (\tau )$
touches the boundary 
$$
\partial B\left( x,\left| f_{ij}(x)\right| \left/
6C_{23}\left( d(x,x_0)+1\right) ^{k_2}\right. \right),
$$
it is easy to see that
\begin{eqnarray}
\frac{\left|
f_{ij}(x)\right| }{6C_{23}\left( d(x,x_0)+1\right) ^{k_2}}
& \leq & \mbox{the length of} \; \gamma _\xi \; \mbox{on} \;
[0,\tau_0] \nonumber\\
& \leq & 2C_{22}\int_0^{\tau _0}\frac{\left( d(x,x_0)+1\right)
^{k_1}}{\left| f_{ij}(x)\right| }d\tau \qquad (\mbox{by} \; (8.7))
\nonumber\\
& = & 2C_{22}\tau _0\frac{\left( d(x,x_0)+1\right) ^{k_1}}{\left|
f_{ij}(x)\right| }\ .\nonumber
\end{eqnarray}
Thus,
\begin{equation}
\label{8.11}
\tau _0\geq \frac{\left| f_{ij}(x)\right| ^2}{
2C_{22}C_{23}6\left( d(x,x_0)+1\right) ^{k_1+k_2}}\ . 
\end{equation}

Note that the integral curve $\gamma_{\xi}$ projects to straight
line passing through $f(x)$ by $\rho$. Thus, when
$(\alpha_i, \alpha_j, \beta_i, \beta_j)$ runs through the
unit sphere in $\C^2$, the collection of integral curves 
$\gamma_{\xi}$ inside  
$$
B\left(x,\left|f_{ij}(x)\right|\left/6C_{23}\left(d(x,x_0)+1
\right) ^{k_2}\right.\right)
$$ 
will project, by $\rho$, onto the Euclidean ball
$$
B_{\C^2}\left( \rho(x),\left| f_{ij}(x)\right|^2 
\left/ 2C_22C_26 \left( d(x,x_0)+1\right) ^{k_1+k_2}\right. \right).
$$
This proves (\ref{8.8}) and hence the lemma.
\hfill$\Box $

%Let us define a map $\Phi $ from the unit complex sphere of $\C^2$ to $%
%\partial B(x,|f_{ij}(x)|/\\6^{k_2}C_{23}\left( d(x,x_0)+1\right) ^{k_2})$ by
%sending $(\alpha _1-\sqrt{-1}\beta _1,\alpha _2-\sqrt{-1}\beta _2)$ to $%
%\gamma _\xi (\tau _0)$. By the smooth dependence of the solution of ODE, $%
%\Phi $ is a smooth map. Since $df_i\wedge df_j$ has no zero in $B\left(
%x,\left| f_{ij}(x)\right| \left/ 6^{k_2}C_{23}\left( d(x,x_0)+1\right)
%^{k_2}\right. \right) $ by (\ref{8.6}), the differential of $\Phi $ is
%nondegenerate. In particular, $\Phi $ is open and is a surjective map. This
%says that any boundary point can be joined with $x$ by some integral curve $%
%\gamma _\xi $ in $B\left( x,\left| f_{ij}(x)\right| \left/
%6^{k_2}C_{23}\left( d(x,x_0)+1\right) ^{k_2}\right. \right) $. Thus by
%combining with (\ref{8.10}) and (\ref{8.11}) it is not hard to see that (\ref
%{8.8}) holds. Hence we have complete the proof of Lemma 8.1.

\vskip 3mm 
Now, we are ready to prove the almost surjectivity of the
holomorphic map $F:M\rightarrow {\C^4}$. For each $1\leq i,j\leq 4$, 
since $f_{ij}$ is a holomorphic function of polynomial growth 
and $R(M)$ is generated by $f_1,...,f_4$, we can write
$$
f_{ij}(x) = H_{ij}\left( f_1(x),f_2(x),f_3(x),f_4(x)\right) \qquad 
\mbox{on} \quad M,
$$
for some rational function $H_{ij}$ on ${\C^4}$. Let $T_0$ be the 
union of the singular set of $Z$ and the zero and pole sets of all $H_{ij},\ 
1\leq i,j\leq 4$. For any $b\in Z\backslash (F(M)\cup T_0)$, there exist 
fixed $\{i,j\}\subset \{1,2,3,4\}$ such that the projection 
$\pi _{ij}:Z\rightarrow {\C^2}$ is nondegenerate at $b$. 
Since $Z$ is algebraic, the number of points contained in 
$\pi _{ij}^{-1}\circ \pi _{ij}(b)$ is less than some fixed integer 
$K$ depending only on $Z$. By interpolation, there is a
polynomial $h_b$ of degree $\leq K$ on ${\C^4}$ such that $h_b(b)=1$, 
and $h_b(w)=0$ for all $w\in (\pi _{ij}^{-1}\circ \pi _{ij}(b))\backslash \{b\}$. 
We now solve on $M\backslash V_0$ the ideal problem with unknown 
holomorphic functions $g_i$ and $g_j$,
\begin{equation}
\label{8.12}
\left( f_i-b_i\right) g_i+\left( f_j-b_j\right) g_j=\left(
h_b\circ F\right) ^4\ , 
\end{equation}
where $b=(b_1,b_2,b_3,b_4).$

Let 
$$
\psi =-n_1\log \delta +n_2\log (1+|f_i|^2+|f_j|^2),
$$ 
where the integers $n_1,\ n_2 > 0$ will be determined later. 
Clearly, $\psi $ is a strictly plurisubharmonic function on $M\backslash V_0$. 
By the estimate of Skoda (c.f. Theorem 1.3 in \cite{Mo3}), given any 
$\alpha >1$, there exists a solution $\{g_i,g_j\}$ to (\ref{8.12}) such that 
\begin{eqnarray}
\label{8.13}
& &\int\limits_{M\backslash V_0}\frac{\left( \left| g_i\right|^2+
\left| g_j\right| ^2\right) e^{-\psi }}{\left( \left| f_i-b_i\right|^2+
\left| f_j-b_j\right| ^2\right) ^{2\alpha }}\rho ^{*}dV_E \nonumber\\
& \leq & C_\alpha\int\limits_{M\backslash V_0}\frac{\left( h_b\circ F\right) ^8
e^{-\psi }}{\left( \left| f_i-b_i\right| ^2+\left| f_j-b_j\right| ^2
\right) ^{2\alpha +1}}\rho^{*}dV_E\ ,
\end{eqnarray}
provided the right hand side is finite. Recall that 
$\rho = \pi_{ij} \circ F$ and here 
$$
\rho ^{*}dV_E=\pm \left( \frac{\sqrt{-1}}2\right) ^2df_i\wedge 
\overline{df}_i\wedge df_j\wedge \overline{df}_j
$$ 
denotes the pull back of the Euclidean volume element of ${\C^4.}$

Let $\{\zeta _1,\zeta _2,\cdots ,\zeta _m\}=\pi _{ij}^{-1}\circ \pi _{ij}(b)$
$(m<K)$ be the preimages of $\pi _{ij}(b)$ with $\zeta _1=b$. And let 
$U_k\ ( 1\leq k\leq m )$ be disjoint small neighborhoods of $\zeta _k\
( 1\leq k\leq m )$. The integral on the right hand side of (\ref{8.13})
can be decomposed into three parts
\begin{eqnarray}
RHS
& = &\left( \int_{F^{-1}(U_1)}+\sum_{k=2}^m\int_{F^{-1}(U_k)}+
\int_{\left( M\backslash V_0\right) \left\backslash 
\cup _{k=1}^mF^{-1}(U_k)\right. }\right) \nonumber\\
&   &\frac{\left( h_b\circ F\right) ^8e^{-\psi }}{\left( 
\left| f_i-b_i\right| ^2+\left| f_j-b_j\right| ^2\right) ^{2\alpha +1}}
\rho ^{*}dV_E \nonumber\\
& = & I_1+I_2+I_3\ .\nonumber
\end{eqnarray}

For $I_1$, since $h_b(b)=1$ and $\delta (x)\leq \left( \left| f_i-b_i\right|
^2+\left| f_j-b_j\right| ^2\right) ^{\frac 12}$, we can choose $n_1\geq
2(2\alpha +1)$ and $U_1$ small enough so that the integral $I_1$ is finite.

For $I_2$, since $h_b(\zeta _k)=0$ for $2\leq k\leq m,$ we can choose $%
\alpha $ such that $2(2\alpha +1)<8$ (e.g. $\alpha =1.4$). Then the
integral $I_2$ is also finite.

For $I_3$, we choose $n_2\geq 10+8K+n_1$, where $h_b$ is of degree $\leq K$.
Then $I_3$ can be estimated as
$$
I_3\leq C_{25}\int_{{\C^2}}\frac 1{\left( 1+|w|^2\right) ^{10}}dV_E
<+\infty\ . 
$$

Hence, we have obtained a solution $\{g_i,g_j\}$ of the ideal problem 
(\ref{8.12}) such that 
\begin{equation}
\label{8.14}
\int\limits_{M\backslash V_0}\frac{\left( \left| g_i\right|
^2+\left| g_j\right| ^2\right) e^{-\psi }}{\left( \left| f_i-b_i\right|
^2+\left| f_j-b_j\right| ^2\right) ^{2\alpha }}\rho ^{*}dV_E<+\infty \ .
\end{equation}
Recall from Lemma 8.1 and (\ref{8.5}), we have
$$
\delta (x)\geq C\left| f_{ij}(x)\right| ^2\left( d(x,x_0)+1\right) ^{-p} 
$$
and%
\begin{eqnarray}
\rho ^{*}dV_E&=&\pm \left( \frac{\sqrt{-1}}2\right) ^2df_i\wedge df_j\wedge \overline{df}_i\wedge \overline{df}_j\nonumber\\
&\geq&\frac{\left| f_{ij}\right| ^2}{\left| v\wedge \overline{v}\right| }
\omega ^2\nonumber\\
&\geq&C_{26}\frac{\left| f_{ij}(x)\right| ^2}{\left( d(x,x_0)+1\right) ^{k_3}}\omega ^2\nonumber
\end{eqnarray}
for some positive constants $C_{26}$ and $k_3$. Substituting
these two inequalities into (\ref{8.14}) we get 
\begin{equation}
\label{8.15}
\int\limits_{M\backslash V_0}\frac{\left( \left| g_i\right|
^2+\left| g_j\right| ^2\right) \left| f_{ij}\right| ^{2+2n_1}}{\left(
d(x,x_0)+1\right) ^{k_4}}\omega ^2<+\infty \ 
\end{equation}
where $k_4$ is some positive constant independent of $b$ and $i,\ j.$ Then,
both $g_if_{ij}^{n_1+1}$ and $g_jf_{ij}^{n_1+1}$ are locally
square integrable. They can thus be extended holomorphically from 
$M\backslash V_0$ to $M$. By the mean value inequality of
subharmonic functions, we deduce also that they are of polynomial
growth with degree bounded by some positive number $k_5$ 
independent of $b$.

Now, recall that $R(M)={\C}(f{_1,}f{_2,}f{_3/}f{_4)}$,
the holomorphic functions $g_if_{ij}^{n_1+1}$ and $g_jf_{ij}^{n_1+1}$ 
are thus rational functions of $f{_1,\ }f{_2,\ }f{_3}$ and $f{_4}$. 
Hence, we can regard the equation (\ref{8.12}) as an equation on the 
variety $Z\subset {\C^4}$, namely
$$
\left( w_i-b_i\right) g_if_{ij}^{n_1+1}+\left( w_j-b_j\right)
g_jf_{ij}^{n_1+1}=H_{ij}^{n_1+1}\cdot h_b^4\ . 
$$
Since $h_b$ is a polynomial with $h_b(b)=1$ and the point $b$ lies outside
of the zero and pole sets of $H_{ij}$, either $g_if_{ij}^{n_1+1}$ or 
$g_jf_{ij}^{n_1+1}$, when regarded as rational function on $\C^4$, 
must have a pole at $b$. Denote this function by $G^0$. 
Thus, $G^0$ is a rational function on $Z$ with $G^0(b)=\infty $ and 
$G^0\circ F$ is a holomorphic function on $M$ with degree $\leq k_5$. 
If $Z\backslash (F(M)\cup T_0 \cup \mbox{pole sets of}\, G^0)$ is empty,
then we are done. Otherwise, pick any 
$b_1\in Z\backslash (F(M)\cup T_0 \cup \mbox{pole sets of}\, G^0)$ and 
repeat the same procedure to obtain a rational function $G^1$ on $Z$ with 
$G^1(b_1)=\infty $ and $G^1\circ F$ a holomorphic function on $M$ with
degree $\leq k_5$. Proceeding this way, we obtain a
sequence of points $\{b,b_1,b_2,\cdots \}$ and rational functions 
$\{G^0,G^1,G^2,\cdots \}$ such that $G^k(b_k)=\infty $ and $G^l$ regular
at $b_k$ for $l<k$. So, $\{G^0,G^1,G^2,\cdots \}$ must be linearly independent
over ${\C}$. Moreover, all of $G^k\circ F$ are holomorphic functions with
degree $\leq k_5$. Hence, by (\ref{8.2}), the above procedure must terminate 
in a finite number of steps. In other words, there exists an algebraic 
subvariety $T$ of $Z$ such that $F(M)\supset Z\backslash T$. 

Moreover, $F$ establishes a quasi embedding from $M$ to a quasi--affine 
algebraic variety. Indeed, let $W=F^{-1}(T)$. By the definition of $T_0$ and 
the construction of $T$, we know that $W \supset V$, where $V$ is the 
union of the branching locus of $F$ and $F^{-1}(\mbox{Sing}(Z))$, and 
$W$ is the zero divisor of finitely many holomorphic functions of 
polynomial growth. Therefore, $F$ maps $M\backslash W$ biholomorphically 
onto $Z\backslash T.$

Finally, to complete the proof of our Main Theorem, we have to show that 
the mapping $F$ can be desingularized by adjoining a finite number of 
holomorphic functions of polynomial growth and taking normalization of
the image. 

We have constructed the mapping $F:M\rightarrow Z$ into an affine algebraic
variety which maps $M\backslash W$ biholomorphically onto $Z\backslash T.$
Now, we use normalization of the affine algebraic variety $Z$ to resolve the
codimension $1$ singularities of $F.$ 
Let $\mbox{Reg}(Z)$ denote the Zariski dense
subset of $Z$ consisting of its regular points. It is well known that the
normalization $\widetilde{Z}$ of $Z$ can be obtained by taking 
$\widetilde{Z}$ to be the closure of the graph of $\{Q_{1,}Q_2,\cdots ,Q_m\}$ 
on $\mbox{Reg}(Z)$
where $Q_i$ is a rational function which is holomorphic (or regular in the 
terminology of algebraic geometry) on $\mbox{Reg}(Z)$. The lifting of 
$F:M\rightarrow Z$ to $\widetilde{F}:M\rightarrow Z$ is then given by 
$\{f{_1,}f{_2,}f{_3,}f{_4,}Q_1\circ F,\cdots ,Q_m\circ F\}$ where, as 
was shown in proposition 8.1 of Mok \cite{Mo1}, for each $i$, $Q_i\circ F$ 
can be holomorphically extended to the whole manifold $M$ as a 
holomorphic function of polynomial growth. 

Write $F_0=F:M\rightarrow Z$ and
denote $\widetilde{F}_0:M\rightarrow \widetilde{Z}$ the normalization of 
$F_0$. For any smooth point $x$ on the subvariety $W$, by using the 
$L^2$ estimates of the $\overline{\partial }$ operator as in Section 7, 
one can find two holomorphic functions $g_x^1,\ g_x^2$ of polynomial growth 
which give local holomorphic coordinates at $x$. Adding $g_x^1,\ g_x^2$ to 
the map $\widetilde{F}_0$, we get a new map 
$F_1=(\widetilde{F}_0,g_x^1,g_x^2):M\rightarrow Z_1\subset {{\C}^{6+m}}$, 
which is nondegenerate at $x$. 
Write the normalization of $F_1$ as 
$\widetilde{F}_1:M\rightarrow \widetilde{Z}_1$ and
continue in this way to get holomorphic mappings $F_i:M\rightarrow Z_i$ and
their normalizations $\widetilde{F}_i:M\rightarrow \widetilde{Z}_i$ such that%
$$
\widetilde{W}_0\stackrel{\supset }{\neq }\widetilde{W}_1\stackrel{\supset }{%
\neq }\cdots \stackrel{\supset }{\neq }\widetilde{W}_i\stackrel{\supset }{
\neq }\cdots,  
$$
where $\widetilde{W}_i$ is the locus of ramification of $\widetilde{F}_i.$

Note that $\widetilde{W}_i$ contains no isolated
point because $\widetilde{Z}_i$ is normal. Moreover, by  
Proposition 7.2, $W$ has only finite number of
irreducible components because $W$ is the zero divisor of finitely
many holomorphic
functions of polynomial growth. This implies that the above 
procedure must terminate in a finite number of steps, say $l$. 
%Note that $\widetilde{V}_l$ contains no isolated
%points because $\widetilde{Z}_l$ is normal. 
%Hence $\widetilde{V}_l=\emptyset $. 
Thus, we get a biholomorphism $\widetilde{F}_l$ from $M$ onto its image 
$\widetilde{F}_l(M)\subset \widetilde{Z}_l$. The argument in our proof 
of the almost surjectivity shows that $\widetilde{F}_l(M)$ can miss at 
most finitely many irreducible subvarieties of $\widetilde{Z}_l$, 
say $\widetilde{T}_1^{(l)},\cdots ,\widetilde{T}_q^{(l)}$. 
If $\widetilde{F}_l(M)\cap \widetilde{T}_i^{(l)}\neq \emptyset $, 
then it must intersect $\widetilde{T}_i^{(l)}$ in a nonempty open set 
because $\widetilde{F}_l$ is open. We arrange $\widetilde{T}_i^{(l)}$ 
so that $\widetilde{F}_l(M)\cap \widetilde{T}_i^{(l)}=\emptyset $ 
for $1\leq i\leq p$ and 
$\widetilde{F}_l(M)\cap \widetilde{T}_i^{(l)}\neq \emptyset $ for 
$p+1\leq i\leq q$. Note that$\widetilde{F}_l(M)$ is a Stein subset 
of $\widetilde{Z}_l$ because $M$ is Stein by Theorem 5.1 and 
$\widetilde{F}_l$ maps $M$ biholomorphically onto its image. 
By Hartog's extension theorem, every holomorphic function on 
$\widetilde{Z}_l\backslash \cup _{1\leq i\leq q}\widetilde{T}_i^{(l)}$
extends to 
$\widetilde{Z}_l\backslash \cup _{1\leq i\leq p}\widetilde{T}_i^{(l)}$. 
Hence, we get a biholomorphic map from $M$ onto a quasi--affine
algebraic variety. Finally, recall that a classical theorem of Ramanujam \cite{R} in
affine algebraic geometry says that an algebraic variety homeomorphic to 
${\R^4}$ is biregular to ${\C^2}$. Combining this result of Ramanujam with
Theorem 5.1, we deduce that $M$ is actually biholomorphic to ${\C^2}$.
Therefore we have completed the proof of the Main Theorem.
%\hfill$\Box $

\end{document}